\documentclass[11pt]{article}

\usepackage[T1]{fontenc}
\usepackage[utf8]{inputenc}
\usepackage{lmodern}
\usepackage{microtype}
\usepackage[a4paper,margin=1.12in]{geometry}
\usepackage{amsmath,amssymb,amsthm,mathtools}
\usepackage{enumitem}
\usepackage[hidelinks]{hyperref}
\usepackage[nameinlink,noabbrev,capitalise]{cleveref}
\usepackage{color}

\allowdisplaybreaks[1]

\newtheorem{theorem}{Theorem}[section]
\newtheorem{proposition}[theorem]{Proposition}
\newtheorem{lemma}[theorem]{Lemma}
\newtheorem{corollary}[theorem]{Corollary}
\theoremstyle{definition}

\theoremstyle{remark}
\newtheorem{remark}[theorem]{Remark}

\newcommand{\Z}{\mathbb Z}
\newcommand{\Q}{\mathbb Q}
\newcommand{\R}{\mathbb R}
\newcommand{\C}{\mathbb C}
\newcommand{\T}{\mathbb T}
\newcommand{\1}{\mathbf 1}
\newcommand{\supp}{\operatorname{supp}}
\newcommand{\Gal}{\operatorname{Gal}}
\newcommand{\lcm}{\operatorname{lcm}}
\newcommand{\diag}{\operatorname{diag}}

\newcommand{\numberthis}[1]{\eqno\hbox{\normalfont
  \refstepcounter{equation}\label{#1}(\theequation)}}

\title{Periodicity of weak tilings of integers}
\author{Jialun Li\thanks{\raggedright School of Mathematical Sciences, Fudan University, Shanghai 200433, China. Email: \texttt{jialunli@fudan.edu.cn}.}
\and
Ruxi Shi\thanks{\raggedright Shanghai Center for Mathematical Sciences, Fudan University, Shanghai 200438, China. Email: \texttt{ruxishi@fudan.edu.cn}. R.S. was supported by the NSFC Nos.~12571198 and 12231013.}}
\date{}

\begin{document}
\maketitle

\begin{abstract}
We prove that every finite set of integers admitting a weak tiling $x$
with $x(0)=1$ also admits a periodic weak tiling with the same value at
the origin.
The periodic tiling is obtained by retaining only the root-of-unity
frequencies of $x$.  The proof combines Galois filtering with positive
semidefinite matrix estimates; an alternative argument uses harmonic
extension on a strip.  We deduce periodicity of normalized
positive-definite weak tilings and recover the rationality theorem for
finite spectra.
An example shows that normalized weak tilings can be nonperiodic,
even though their sets of points attaining the value $1$ are always periodic.
\end{abstract}

\section{Introduction}
\label{sec:introduction}

Every translational tiling of $\Z$ by a finite set is periodic, by the
classical theorem of Newman~\cite{Newman}.  In higher dimensions, the
\emph{periodic tiling conjecture} asserts that every finite set
$A\subset\Z^d$ that tiles $\Z^d$ by translations also admits a periodic
tiling, with a translation set invariant under a finite-index subgroup
of $\Z^d$.  Bhattacharya proved the conjecture for $d=2$~\cite{Bhattacharya},
whereas Greenfeld and Tao disproved it for sufficiently large $d$ by
constructing a finite set that tiles $\Z^d$ but admits no periodic
tiling~\cite{GreenfeldTao}.

Replacing the indicator of the translation set by a nonnegative weight
leads to the notion of a
\emph{weak tiling}: for a finite nonempty set $A\subset\Z$, this is a
function $x:\Z\to[0,\infty)$ satisfying the convolution equation
$$
 (\1_A*x)(n)=\sum_{a\in A}x(n-a)=1\qquad\text{for all }n\in\Z.
$$
They need not be periodic.

Weak tiling provides a useful connection between translational tilings
and spectral sets.  Lev and Matolcsi showed that a bounded spectral
set in $\R^d$ weakly tiles its complement, a central step in their proof of
Fuglede's conjecture for convex bodies~\cite{LevMatolcsi}.
Further properties of this condition were studied by Kolountzakis,
Lev, and Matolcsi~\cite{KLM}.  In finite abelian groups, positive-definite
weak tilings were studied in~\cite{KMMS}, and Kiss, Londner, Matolcsi,
and Somlai constructed a set that weakly tiles its complement but is
neither spectral nor a translational tile~\cite{KLMS}. For broader
background on spectral sets and weak tiling, see~\cite{KolountzakisSurvey}.
 Our question
concerns the periodicity of weak tilings normalized by $x(0)=1$.  Equivalently,
in the terminology recalled below, we ask whether a finite set that
weakly tiles its complement always admits a periodic weak tiling of
its complement.

\begin{theorem}\label{thm1}
Let $A\subset\Z$ be finite.  If $A$ admits a weak tiling
$x$ with $x(0)=1$, then it also admits a periodic weak tiling $y$
with $y(0)=1$.
\end{theorem}
\subsection*{Set-up and precise statements}

We now take a closer look at weak tilings, starting with basic definitions. For a finite nonempty set $A\subset\Z$, define
$$
 X_A=\left\{x:\Z\to[0,1]:
       \sum_{a\in A}x(n-a)=1\quad(n\in\Z)\right\}.
$$
The set $X_A$ consists of the weak tilings of $\Z$ by $A$,
that is, the solutions of $\1_A*x\equiv1$.  A nonnegative solution
of this equation automatically satisfies $x\le1$, so nothing is lost
by restricting the values to $[0,1]$. The space $X_A$ is compact in the product topology inherited from
$[0,1]^\Z$ and invariant under the shift $(\sigma(x))(n)=x(n+1)$.
Thus $(X_A,\sigma)$ is a topological dynamical system. A point $x\in X_A$ is \emph{periodic} if there is an integer $N\ge1$ such that $\sigma^N(x)=x$, i.e. $x(n+N)=x(n)$ for every $n\in\Z$. The constant sequence $x\equiv1/|A|$ always belongs to $X_A$;
in particular, the system $(X_A,\sigma)$ is nonempty and contains a
periodic point. 

We are interested in the normalized subset
$$
 X_A^{(1)}=\{x\in X_A:x(0)=1\}.
$$
This is a closed subset of $X_A$, but it may be empty.
Each element of $X_A^{(1)}$ corresponds to a weak tiling of the
complement of $A$ in the terminology of~\cite{LevMatolcsi,KLMS}.
Indeed, for $x\in X_A^{(1)}$ the weight
$\mu=x-\delta_0$, where $\delta_0=\1_{\{0\}}$, satisfies
$$
 \mu\ge0,\qquad \1_A*\mu=\1_{\Z\setminus A}.
$$
Conversely, if a nonnegative $\mu$ satisfies this identity, evaluating
it at any $a\in A$ gives $\mu(0)=0$, and $x=\mu+\delta_0$ belongs to
$X_A^{(1)}$.  Thus $X_A^{(1)}\ne\varnothing$ if and only if $A$
weakly tiles its complement.


\begin{theorem}\label{thm:main}
Let $A\subset\Z$ be finite and nonempty.  If $X_A^{(1)}$ is nonempty,
then it contains a periodic point.
\end{theorem}

\cref{thm:main}  is equivalent to \cref{thm1} and both are consequences of a
more precise construction: the periodic weak tiling is obtained
canonically from the original weight by retaining only the
root-of-unity frequencies in its exponential expansion. The next
theorem makes this precise.  Write
$\T=\{z\in\C:|z|=1\}$ for the unit circle.

\begin{theorem}\label{thm:expansion}
Let $A\subset\Z$ be finite and nonempty.  Every $x\in X_A$ has a unique
representation
$$
 x(n)=\frac1{|A|}+\sum_{\lambda\in\Lambda_x}c_\lambda\lambda^n,
 \qquad n\in\Z,
\numberthis{eq:expansion-intro}
$$
where $\Lambda_x\subset\T\setminus\{1\}$ is finite and the coefficients
$c_\lambda\in\C$ are nonzero. Let
$$
 y(n)=\frac1{|A|}+\sum_{\substack{\lambda\in\Lambda_x\\
                    \lambda\text{ is a root of unity}}}
                    c_\lambda\lambda^n.
\numberthis{eq:ru-projection-intro}
$$
Then the following properties hold.
\begin{itemize}
    \item [(1)] $y\in X_A$;
    \item [(2)] If $x(0)=1$, then $y\in X_A^{(1)}$ and is periodic;
    \item [(3)] The set
$
 \{n\in\Z:x(n)=1\}$ is equal to $\{n\in\Z:y(n)=1\}
$, which is periodic.
\end{itemize}
\end{theorem}


As we were completing this paper, we learned that Chun-Kit Lai had asked
whether every weak tiling $x$ with $x(0)=1$ must be periodic. He ruled
out nonperiodic solutions in several special cases~\cite{Lai}. However,
the origial
weight \eqref{eq:expansion-intro} need not be periodic: \cref{prop:nonperiodic-example} gives a
nonperiodic $x\in X_A^{(1)}$ with $\gcd(A-A)=1$.  Also, the periodic
weight in~\eqref{eq:ru-projection-intro} need not be $\{0,1\}$-valued,
so the conclusion does not assert a proper translational tiling.

\paragraph*{Idea of proof} The proof adapts the hyperbolic cancellation and positive-semidefinite
comparison in Fu and Song's modulus-rigidity argument for finite
Hadamard pairs~\cite[Theorem~1.1, (4.11)--(4.12)]{FuSong}.
Their argument uses the special inverse-orthogonality structure of
finite Hadamard pairs.  A general $x\in X_A^{(1)}$ need not be a
spectral weight, and the coefficients in its finite exponential
expansion need not be nonnegative.  We therefore establish a comparison
estimate for general nonnegative functions satisfying the weak-tiling
identities, using positivity of their Fourier-coefficient matrices
in place of nonnegativity of individual coefficients.

 The finite expansion~\eqref{eq:expansion-intro} follows from the bounded
linear recurrence satisfied by $x-1/|A|$, as in the work of Akhmedov
and Fulghesu on arithmetic sets~\cite{AkhmedovFulghesu}.  In fact, the frequencies in~\eqref{eq:expansion-intro} are roots of the
polynomial $P_A(z)=\sum_{a\in A}z^a$, and hence algebraic integers.
Let $K/\Q$ be a finite Galois extension containing them.
For each $\sigma\in\Gal(K/\Q)$, define
$$
 y_\sigma(n)=\frac1{|A|}
 +\sum_{\substack{\lambda\in\Lambda_x\\|\sigma(\lambda)|=1}}
 c_\lambda\lambda^n.
$$
We call this operation a \emph{Galois filter}: it retains the term
$c_\lambda\lambda^n$ precisely when $|\sigma(\lambda)|=1$, keeping
its original coefficient.
Using a compact-group representation and the matrix comparison,
we show that $y_\sigma\in X_A$ and
$$
 0\le c_A\bigl(1-y_\sigma(n)\bigr)\le1-x(n)
 \qquad\text{for all }n\in\Z,
$$
for a constant $c_A>0$ depending only on $A$;
see~\cref{lem:filter,prop:one-filter}.
Thus $x(n)=1$ implies $y_\sigma(n)=1$, so each filter preserves
the normalization at the origin.
An alternative proof of this estimate uses harmonic extension on a
strip and the minimum principle; see~\cref{sec:filter-harmonic}.

Applying these filters successively for all $\sigma\in\Gal(K/\Q)$
leaves precisely the frequencies whose algebraic conjugates all have
modulus $1$.   By Kronecker's theorem, these algebraic integers
are roots of unity.  The resulting sequence is
therefore the periodic weak tiling in~\eqref{eq:ru-projection-intro}.
An additional consequence is that a normalized
positive-definite weight is itself periodic.  Applied to spectral
weights, this recovers the finite rational-spectrum theorem of
Fu and Song~\cite[Theorem~1.2]{FuSong}; see~\cref{cor:rational-spectrum}.

\paragraph*{Structure of the paper}We recall the matrix positivity tools in~\cref{sec:positivity} and
establish the recurrence and compact-group representation
in~\cref{sec:preliminaries}.  The one-step Galois filter and its
pointwise estimate are proved in~\cref{sec:filter}, and the iteration
is completed in~\cref{sec:main-proof}.  Consequences and applications
are given in~\cref{sec:consequences,sec:applications}, followed by
the nonperiodic example in~\cref{sec:example}.  In~\cref{sec:dynamics}
we prove that the shift on $X_A$ has zero topological entropy and
describe its orbit closures, its invariant measures, and the set $X_A^{(1)}$.

\section{Positive semidefinite matrices and the hyperbolic kernel}
\label{sec:positivity}

For a Hermitian matrix $R$, the notation
$$
 R\succeq0
$$
means that $R$ is positive semidefinite, that is, $v^*Rv\ge0$ for every
vector $v$.  For a finite index set $F$, we write $I_F$ for the identity
matrix indexed by $F$.

We use the Schur product theorem~\cite{Schur}.

\begin{theorem}\label{thm:schur-product}
Let $R$ and $S$ be Hermitian matrices of the same size.  If $R\succeq0$ and $S\succeq0$, then their entrywise product $R\circ S$ satisfies $R\circ S\succeq0$.
\end{theorem}

The following classical hyperbolic kernel is used in the modulus-rigidity
argument of Fu and Song~\cite[Lemma~4.1]{FuSong}.  Its positive
definiteness and Fourier representation appear, after rescaling,
in Bhatia and Parthasarathy~\cite[Theorem~3.2]{BhatiaParthasarathy}.
For $0<s<M$, define
$$
 k_{s,M}(t)=
 \begin{cases}
 \dfrac{\sinh(st)}{\sinh(Mt)},&t\ne0,\\[6pt]
 \dfrac{s}{M},&t=0.
 \end{cases}
\numberthis{eq:kernel}
$$
It admits the positive Fourier representation
$$
 \begin{aligned}
 k_{s,M}(t)&=\int_{\R}e^{i\xi t}\,d\mu_{s,M}(\xi),\\
 d\mu_{s,M}(\xi)&=
 \frac{\sin(\pi s/M)}{2M}
 \frac{d\xi}{\cosh(\pi\xi/M)+\cos(\pi s/M)}.
 \end{aligned}
\numberthis{eq:positive-measure}
$$

The density in~\eqref{eq:positive-measure} is a rescaled Poisson
kernel for one boundary component of a strip; see
Widder~\cite{Widder} and~\cite[Section~II.A]{ANNRW}.
This gives a harmonic-measure interpretation of its positivity.

\begin{lemma}\label{lem:hyperbolic}
For $0<s<M$, the function $k_{s,M}$ defined in~\eqref{eq:kernel} is continuous and positive definite on $\R$.
\end{lemma}

\begin{proof}
For $0<s<M$, the density in~\eqref{eq:positive-measure} is positive and decays exponentially at
infinity, so the measure is positive and finite.  Its total mass is
$k_{s,M}(0)=s/M$.  In particular, for real numbers $t_\alpha$ and
complex numbers $c_\alpha$,
$$
 \sum_{\alpha,\beta}\overline{c_\alpha}c_\beta
 k_{s,M}(t_\alpha-t_\beta)
 =\int_{\R}
   \left|\sum_\alpha c_\alpha e^{-i\xi t_\alpha}\right|^2
   d\mu_{s,M}(\xi)\ge0.
$$
This proves positive definiteness, and continuity follows from the
integral representation.
\end{proof}

\section{The recurrence and the compact-group representation}
\label{sec:preliminaries}

Translating $A$ leaves $X_A$ unchanged.  For the proof we therefore normalize $A$ by
$$
 0=\min A,\qquad M=\max A.
\numberthis{eq:normalization}
$$
The case $|A|=1$ gives $x\equiv1$.  Unless stated otherwise, we assume $|A|\ge2$ and $M\ge1$ below.

Define the mask polynomial
$$
 P_A(z)=\sum_{a\in A}z^a
$$
and the characteristic polynomial
$$
 Q_A(z)=z^M P_A(z^{-1})=\sum_{a\in A}z^{M-a}.
\numberthis{eq:characteristic}
$$
Both the leading coefficient and the constant coefficient of $Q_A$ are equal to $1$.

\subsection{Bounded recurrence solutions}

\begin{lemma}\label{lem:recurrence}
Let $u:\Z\to\C$ satisfy
$$
 \sum_{a\in A}u(n-a)=0\qquad(n\in\Z)
\numberthis{eq:homogeneous-recurrence}
$$
and suppose $u$ is bounded.  Then the following hold:
\begin{enumerate}[label=(\roman*),leftmargin=2.2em]
\item There are a unique finite set $\Lambda$ and unique nonzero coefficients $c_\lambda\in\C$ such that
$$
 u(n)=\sum_{\lambda\in\Lambda}c_\lambda\lambda^n,
\numberthis{eq:homogeneous-expansion}
$$
where $\Lambda$ is a subset of
$$
 \{\lambda\in\C:Q_A(\lambda)=0,\ |\lambda|=1\}.
$$
\item For every $\lambda\in\Lambda$ one has
$$
 P_A(\lambda)=P_A(\lambda^{-1})=0.
\numberthis{eq:two-mask-zeros}
$$
In particular, every $\lambda\in\Lambda$ is an algebraic unit, meaning that both $\lambda$ and $\lambda^{-1}$ are algebraic integers.
\item Let $\lambda\in\C$ satisfy $P_A(\lambda)=0$ and $|\lambda|=1$, and let
$p\in\Z[z]$ be its monic minimal polynomial of degree $m$.
Then $p$ divides both $P_A$ and $Q_A$, and
$$
 p(z)=z^m p(z^{-1}).
$$
If $\lambda\ne-1$, then $m$ is even.  In particular, if $P_A$ is
irreducible and has a root $\lambda$ with $|\lambda|=1$, then $P_A=Q_A$ and its degree
is even unless $P_A(z)=z+1$.
\end{enumerate}
\end{lemma}

\begin{proof}
\textup{(i)} Set
$$
 U_n=(u(n),u(n-1),\ldots,u(n-M+1))^{\mathsf T}\in\C^M.
$$
Because $0,M\in A$, equation~\eqref{eq:homogeneous-recurrence} determines both the next and the preceding value of the sequence.  Hence there is an invertible companion matrix $C_A$ with characteristic polynomial $Q_A$ such that
$$
 U_n=C_A^nU_0\qquad(n\in\Z).
$$
Put $C_A$ into Jordan normal form.  A component belonging to an eigenvalue $\lambda$ with $|\lambda|>1$ is unbounded as $n\to+\infty$, while a component with $|\lambda|<1$ is unbounded as $n\to-\infty$.  For $|\lambda|=1$, every nontrivial generalized-eigenvector component grows polynomially in $|n|$.  Since $(U_n)_{n\in\Z}$ is bounded, only ordinary eigenvectors with $|\lambda|=1$ can occur.  This gives~\eqref{eq:homogeneous-expansion}.  Distinct exponential sequences $n\mapsto\lambda^n$ are linearly independent, so the coefficients are unique.

\textup{(ii)} Let $\lambda\in\Lambda$.  By \textup{(i)}, $Q_A(\lambda)=0$ and $|\lambda|=1$.  Equation~\eqref{eq:characteristic} therefore gives $P_A(\lambda^{-1})=0$.  Since $P_A$ has real coefficients and $\overline{\lambda^{-1}}=\lambda$, complex conjugation yields $P_A(\lambda)=0$.  This proves~\eqref{eq:two-mask-zeros}.  In particular, both $\lambda$ and $\lambda^{-1}$ are roots of the monic polynomial $P_A\in\Z[z]$.  Hence both are algebraic integers, so $\lambda$ is an algebraic unit.

\textup{(iii)} Let $\lambda\in\C$ satisfy $P_A(\lambda)=0$ and $|\lambda|=1$.
Since $P_A$ has real coefficients, both $\lambda$ and $\lambda^{-1}$
are zeros of $P_A$.  Thus $Q_A(\lambda)=0$, and minimality gives
$p\mid P_A$ and $p\mid Q_A$.
Likewise $\lambda^{-1}=\overline\lambda$ is a zero of $p$, so
$\lambda$ is a zero of the reciprocal polynomial
$p^*(z):=z^m p(z^{-1})$.  These polynomials have the same degree,
and irreducibility gives $p^*=c p$ for some nonzero $c\in\Q$.
Since $P_A(1)=|A|\ne0$, we have $p(1)\ne0$; evaluating at $1$
yields $c=1$.  If $\lambda\ne-1$, irreducibility also gives
$p(-1)\ne0$.  Evaluating $p=p^*$ at $-1$ then shows $(-1)^m=1$.
\end{proof}

For $x\in X_A$, the sequence
$$
 u(n)=x(n)-\frac1{|A|}
$$
satisfies~\eqref{eq:homogeneous-recurrence}.  Since $0\le x\le1$, it is bounded.  Also $Q_A(1)=P_A(1)=|A|\ne0$, so the constant frequency does not occur in the homogeneous part.  Thus \cref{lem:recurrence} yields the expansion~\eqref{eq:expansion-intro}.

Since $x$ is real-valued, uniqueness also gives
$$
 \lambda^{-1}\in\Lambda_x,\qquad
 c_{\lambda^{-1}}=\overline{c_\lambda}\qquad(\lambda\in\Lambda_x).
$$
\subsection{A compact group determined by the frequencies}

Fix $x\in X_A$, and let $\Gamma=\langle\Lambda_x\rangle\le\T$ be the
multiplicative group generated by its frequencies, equipped with the
discrete topology.  Define the finitely supported coefficient function
$$
 r_x(\gamma)=
 \begin{cases}
  1/|A|,&\gamma=1,\\
  c_\gamma,&\gamma\in\Lambda_x,\\
  0,&\text{otherwise}.
 \end{cases}
$$
The conjugacy relation above gives
$r_x(\gamma^{-1})=\overline{r_x(\gamma)}$.

Let $G=\widehat\Gamma$ be the compact Pontryagin dual of $\Gamma$,
with normalized Haar measure $m_G$; see~\cite{Rudin}.
Its elements are homomorphisms $z:\Gamma\to\T$, and multiplication
is pointwise.  The identity character $e_G$ has $e_G(\gamma)=1$.
The coefficients determine the continuous function
$$
 f_x(z)=\sum_{\gamma\in\Gamma}r_x(\gamma)z(\gamma),\qquad z\in G.
$$
The inclusion of $\Gamma$ in $\T$ is itself a character:
$$
 g(\gamma)=\gamma\qquad(\gamma\in\Gamma).
$$
Consequently,
$$
 x(n)=f_x(g^n)\qquad(n\in\Z).
\numberthis{eq:orbit-values}
$$
This representation allows us to use positivity on a compact group
while retaining all multiplicative relations among the frequencies.

\begin{proposition}\label{lem:compact-dual}
The following hold:
\begin{enumerate}[label=(\roman*),leftmargin=2.2em]
\item The subgroup $\{g^n:n\in\Z\}$ is dense in $G$.
\item We have
$$
 0\le f_x(z)\le1\quad(z\in G),\qquad f_x(e_G)=x(0).
$$
\item For every finite $F\subset\Gamma$, the matrix
$$
 R_F=\bigl(r_x(\alpha\beta^{-1})\bigr)_{\alpha,\beta\in F}
$$
is Hermitian and satisfies $0\preceq R_F\preceq I_F$.
\end{enumerate}
\end{proposition}

\begin{proof}
\textup{(i)} If $\gamma\in\Gamma$ annihilates the closure of
$\{g^n:n\in\Z\}$, then $\gamma^n=g^n(\gamma)=1$ for all $n\in\Z$.
Taking $n=1$ gives $\gamma=1$.  The annihilator is therefore trivial,
and Pontryagin duality gives density.

\textup{(ii)} Use~\eqref{eq:orbit-values}, continuity, and density
to extend $0\le x\le1$ to $G$.  At the identity,
$f_x(e_G)=\sum_\gamma r_x(\gamma)=x(0)$.

\textup{(iii)} Character orthogonality gives
$$
 r_x(\delta)=\int_G f_x(z)\overline{z(\delta)}\,dm_G(z).
$$
Thus, for $c=(c_\alpha)_{\alpha\in F}$,
$$
\begin{aligned}
 0\le c^*R_Fc
 &=\int_G f_x(z)
   \left|\sum_{\alpha\in F}c_\alpha z(\alpha)\right|^2dm_G(z)\\
 &\le\int_G
   \left|\sum_{\alpha\in F}c_\alpha z(\alpha)\right|^2dm_G(z)
   =\sum_{\alpha\in F}|c_\alpha|^2.
\end{aligned}
$$
This proves $0\preceq R_F\preceq I_F$.
\end{proof}

To describe the multiplicative relations explicitly, enumerate
$\Lambda_x=\{\lambda_1,\ldots,\lambda_s\}$.  Evaluation on these
generators identifies $G$ with
$$
 \left\{z\in\T^s:
 \prod_{j=1}^s z_j^{m_j}=1\ \text{whenever}
 \prod_{j=1}^s\lambda_j^{m_j}=1, m\in\Z^s\right\}.
$$
In particular, a phase assignment belongs to $G$ precisely when it
respects every relation among the frequencies.  Since $\Gamma$ is
finitely generated and its torsion subgroup is finite cyclic,
$$
 \Gamma\cong(\Z/q\Z)\oplus\Z^\nu,\qquad
 G\cong\mu_q\times\T^\nu
$$
for some $q\ge1$ and $\nu\ge0$, where $\mu_q$ denotes the group of
$q$th roots of unity.  We will use this structure to approximate
Fourier values by finite matrix quadratic forms.

\section{The hyperbolic kernel and a Galois filter}
\label{sec:filter}

\subsection{Definition of the Galois filter}
\label{sec:filter-setting}

Keep $\Gamma$, $G$, and $g$ from the compactification above, and let
$K$ be the splitting field of $P_A$ over $\Q$; we use the standard
properties of finite Galois extensions as in~\cite{Lang}.  By
\cref{lem:recurrence}, all generators of $\Gamma$ lie in $K$, so
$\Gamma\le K^\times\cap\T$.

For a single filtering step, let $r:\Gamma\to\C$ be finitely supported and satisfy
$$
 r(1)=\frac1{|A|},\qquad
 f(z):=\sum_{\gamma\in\Gamma}r(\gamma)z(\gamma)\in[0,1]\quad(z\in G),
\numberthis{eq:filter-assumptions}
$$
and assume that
$$
 P_A(\gamma)=P_A(\gamma^{-1})=0
 \qquad(\gamma\in\supp r\setminus\{1\}).
\numberthis{eq:root-hypothesis}
$$
These hypotheses hold for $r=r_x$ by
\cref{lem:recurrence,lem:compact-dual}.  We formulate the filtering
step without the condition $f(e_G)=1$, so that it can be iterated on
the same group $\Gamma$.  The anchor condition will enter only after
the pointwise comparison has been proved.

Set
$$
 c_A:=\frac{M-\max(A\setminus\{M\})}{M}\in(0,1].
$$
Thus $1-a/M\ge c_A$ for every $a\in A\setminus\{M\}$.

Fix $\sigma\in\Gal(K/\Q)$.  Define homomorphisms
$$
 \ell(\gamma)=\log|\sigma(\gamma)|\in\R,\qquad
 \omega(\gamma)=\frac{\sigma(\gamma)}{|\sigma(\gamma)|}\in\T
 \qquad(\gamma\in\Gamma).
\numberthis{eq:ell-omega}
$$
The homomorphism $\omega:\Gamma\to\T$ defined in~\eqref{eq:ell-omega}
is a continuous character because $\Gamma$ has the discrete topology.
We therefore regard this map as an element of $G=\widehat\Gamma$, using
the same symbol $\omega$.  In particular, for $z\in G$ and $a\in\Z$,
the product $z\omega^a$ denotes the character
$\gamma\mapsto z(\gamma)\omega(\gamma)^a$.
Both $\ell$ and $\omega$ respect every multiplicative
relation in $\Gamma$.
Let
$$
 \Gamma_0=\ker\ell=\{\gamma\in\Gamma:|\sigma(\gamma)|=1\}.
$$
Define the \emph{filtered coefficient function} and its Fourier sum by
$$
 \begin{aligned}
 r_\sigma(\gamma)&=r(\gamma)\1_{\Gamma_0}(\gamma),\\
 f_\sigma(z)&=\sum_{\gamma\in\Gamma}r_\sigma(\gamma)z(\gamma).
 \end{aligned}
\numberthis{eq:filtered-coefficients}
$$
We call the map $r\mapsto r_\sigma$ the \emph{Galois filter}
associated with $\sigma$.  It retains the coefficient $r(\gamma)$
when $|\sigma(\gamma)|=1$ and sets it to zero otherwise.
For $r=r_x$, this is precisely the Galois filter introduced
in~\cref{sec:introduction}, since
$$
 f_\sigma(g^n)=\frac1{|A|}
 +\sum_{\substack{\lambda\in\Lambda_x\\|\sigma(\lambda)|=1}}
 c_\lambda\lambda^n=y_\sigma(n)
 \qquad\text{for all }n\in\Z.
$$
The formulation for general $r$ allows us to apply the filters
successively on the same frequency group $\Gamma$.

Let
$$
 H=\Gamma_0^\perp
 =\{h\in G:h(\gamma)=1\ \text{for all }\gamma\in\Gamma_0\}.
$$
Let $m_H$ be the normalized Haar measure on $H$.  Since $\Gamma$ is
discrete, $H^\perp=\Gamma_0$.  Character orthogonality gives
$\int_H h(\gamma)\,dm_H(h)=\1_{\Gamma_0}(\gamma)$, and hence
$$
 f_{\sigma}(z)=\int_H f(zh)\,dm_H(h).
\numberthis{eq:conditional-expectation}
$$
Thus restricting the Fourier coefficients to $\Gamma_0$ is equivalent to averaging over its annihilator $H$.
The pointwise estimate in \cref{lem:filter} will show that the filter preserves every point where $f=1$.

\begin{lemma}
\label{lem:filtered-partitions}
The following hold:
\begin{enumerate}[label=(\roman*),leftmargin=2.2em]
\item The filtered function satisfies
$$
 0\le f_{\sigma}(z)\le1\qquad(z\in G).
\numberthis{eq:filtered-range}
$$
\item For every $z\in G$,
$$
 \sum_{a\in A}f_{\sigma}(z\omega^a)=1,
 \qquad
 \sum_{a\in A}f_{\sigma}(zg^{-a})=1.
\numberthis{eq:filtered-partitions}
$$
\end{enumerate}
These assertions do not require $f(e_G)=1$.
\end{lemma}

\begin{proof}
\textup{(i)} By~\eqref{eq:conditional-expectation}, $f_{\sigma}(z)$ is the
average of $f(zh)$ with respect to the probability measure $m_H$.
Since $0\le f\le1$ by~\eqref{eq:filter-assumptions}, this gives
\eqref{eq:filtered-range}.

\textup{(ii)} For $\gamma\in\supp r_{\sigma}\setminus\{1\}$, we have
$\omega(\gamma)=\sigma(\gamma)$, and hence
$$
 P_A(\omega(\gamma))=P_A(\sigma(\gamma))
 =\sigma(P_A(\gamma))=0,
 \qquad P_A(\gamma^{-1})=0.
$$
Since $r_{\sigma}(1)=1/|A|$ and $P_A(1)=|A|$, it follows that
$$
 r_{\sigma}(\gamma)P_A(\omega(\gamma))
 =r_{\sigma}(\gamma)P_A(\gamma^{-1})
 =\1_{\{1\}}(\gamma)
 \qquad(\gamma\in\Gamma).
\numberthis{eq:coefficient-partition}
$$
For every $z\in G$, finite support and~\eqref{eq:coefficient-partition} give
$$
\begin{aligned}
 \sum_{a\in A}f_{\sigma}(z\omega^a)
 &=\sum_{\gamma\in\Gamma}r_{\sigma}(\gamma)z(\gamma)
   \sum_{a\in A}\omega(\gamma)^a\\
 &=\sum_{\gamma\in\Gamma}r_{\sigma}(\gamma)z(\gamma)
   P_A(\omega(\gamma))=z(1)=1,\\[4pt]
 \sum_{a\in A}f_{\sigma}(zg^{-a})
 &=\sum_{\gamma\in\Gamma}r_{\sigma}(\gamma)z(\gamma)
   \sum_{a\in A}\gamma^{-a}\\
 &=\sum_{\gamma\in\Gamma}r_{\sigma}(\gamma)z(\gamma)
   P_A(\gamma^{-1})=z(1)=1.
\end{aligned}
$$
\end{proof}

\subsection{Hyperbolic cancellation and matrix identities}
\label{sec:filter-master}

We begin with the idea of the proof in~\cref{sec:filter-master,sec:filter-deficit}.
Our goal is to show that $f(z_0)=1$ implies $f_\sigma(z_0)=1$.
In this subsection, we combine the root identities
$P_A(\sigma(\gamma))=P_A(\sigma(\gamma)^{-1})=0$ for each nonconstant
frequency $\gamma$ to obtain coefficient relations involving the
hyperbolic kernels $k_{a,M}$.  The positive definiteness of these
kernels and the Schur product theorem then yield a matrix identity
with positive semidefinite summands.
In~\cref{sec:filter-deficit}, this identity gives the comparison
$I_F-R\succeq c_A(I_F-R_\sigma)$ on every finite set $F\subset\Gamma$.
Testing this comparison on character vectors along a F\o lner sequence
and passing to the limit yields $c_A(1-f_\sigma)\le1-f$.
Since $c_A>0$ and $f_\sigma\le1$, this inequality implies
$f_\sigma(z_0)=1$ whenever $f(z_0)=1$, as required. 

\begin{lemma}
\label{lem:coefficient-balance}
For every $\gamma\in\supp r$ with $\ell(\gamma)\ne0$,
$$
 \omega(\gamma)^M+
 \sum_{a\in A\setminus\{0,M\}}
 \frac{\sinh(a\ell(\gamma))}{\sinh(M\ell(\gamma))}
 \omega(\gamma)^a=0.
\numberthis{eq:hyperbolic-cancellation}
$$
In particular, for every $\gamma\in\Gamma$,
$$
 r(\gamma)\left[
  \omega(\gamma)^M+
  \sum_{a\in A\setminus\{0,M\}}
  k_{a,M}(\ell(\gamma))\omega(\gamma)^a
  \right]=
 r_{\sigma}(\gamma)\left[
  \omega(\gamma)^M+
  \sum_{a\in A\setminus\{0,M\}}\frac aM\omega(\gamma)^a
  \right].
\numberthis{eq:coefficient-balance}
$$
Here $k_{a,M}$ is defined by~\eqref{eq:kernel}, with $k_{a,M}(0)=a/M$.
\end{lemma}

\begin{proof}
For $\gamma\in\supp r$ with $t=\ell(\gamma)\ne0$ and
$u=\omega(\gamma)$, we have $\gamma\ne1$, $|u|=1$, and
$\sigma(\gamma)=e^tu$.
Since $P_A\in\Z[z]$, $\sigma(\gamma^{-1})=e^{-t}u^{-1}$, and
$\overline u=u^{-1}$, \eqref{eq:root-hypothesis} gives
$$
\begin{aligned}
 P_A(e^tu)&=\sigma(P_A(\gamma))=0,\\
 P_A(e^{-t}u)
 &=\overline{P_A(e^{-t}u^{-1})}
 =\overline{\sigma(P_A(\gamma^{-1}))}=0.
\end{aligned}
$$
Consequently,
$$
\begin{aligned}
 0&=\frac{P_A(e^tu)-P_A(e^{-t}u)}2=\sum_{a\in A}\frac{e^{at}-e^{-at}}2\,u^a\\
  &=\sinh(Mt)u^M+
    \sum_{a\in A\setminus\{0,M\}}\sinh(at)u^a.
\end{aligned}
$$
Dividing by $\sinh(Mt)\ne0$ proves~\eqref{eq:hyperbolic-cancellation}.

To prove~\eqref{eq:coefficient-balance}, fix $\gamma\in\Gamma$.
If $r(\gamma)=0$, then $r_{\sigma}(\gamma)=0$, so both sides vanish.
If $\ell(\gamma)=0$, then
$$
 r_{\sigma}(\gamma)=r(\gamma),\qquad
 \lim_{s\to0}\frac{\sinh(as)}{\sinh(Ms)}=\frac aM
 \quad(a\in A\setminus\{0,M\}),
$$
and the two sides coincide.  Otherwise, $r(\gamma)\ne0$ and
$\ell(\gamma)\ne0$; by~\eqref{eq:hyperbolic-cancellation},
$$
 r(\gamma)\left[\omega(\gamma)^M+
 \sum_{a\in A\setminus\{0,M\}}
 \frac{\sinh(a\ell(\gamma))}{\sinh(M\ell(\gamma))}
 \omega(\gamma)^a\right]=0,
 \qquad r_{\sigma}(\gamma)=0.
$$
Thus~\eqref{eq:coefficient-balance} holds in every case.
\end{proof}

We now express the coefficient identities~\eqref{eq:coefficient-partition}
and~\eqref{eq:coefficient-balance} as finite matrix identities.
Let $F\subset\Gamma$ be finite and nonempty.  All matrices below are indexed by $F\times F$.  Define
$$
\begin{aligned}
 R_{\alpha\beta}&=r(\alpha\beta^{-1}),
 & (R_{\sigma})_{\alpha\beta}&=r_{\sigma}(\alpha\beta^{-1}),
 \\
 (W_a)_{\alpha\beta}&=\omega(\alpha\beta^{-1})^a
 \quad(a\in A),
 & (K_a)_{\alpha\beta}&=k_{a,M}(\ell(\alpha\beta^{-1}))
 \quad(a\in A\setminus\{0,M\}).
\end{aligned}
\numberthis{eq:matrix-WK}
$$
In particular, $R_\sigma\circ W_0=R_\sigma$.

The Fourier-coefficient argument in \cref{lem:compact-dual} gives $R\succeq0$ and $R_{\sigma}\succeq0$ from the nonnegativity of $f$ and $f_{\sigma}$.  Each $W_a$ is the rank-one positive matrix
$q_aq_a^*$, where $q_a(\alpha)=\omega(\alpha)^a$.  Since
$\ell(\alpha\beta^{-1})=\ell(\alpha)-\ell(\beta)$, the positivity of
$K_a$ follows from \eqref{eq:positive-measure}.  By \cref{thm:schur-product}, the matrices $R\circ W_M$, $R_{\sigma}\circ W_a$, and $K_a\circ R\circ W_a$ are therefore positive semidefinite whenever the corresponding factors are defined.

\begin{lemma}
\label{lem:matrix-partition}
For every finite nonempty $F\subset\Gamma$,
$$
 \sum_{a\in A}R_{\sigma}\circ W_a=I_F.
\numberthis{eq:matrix-partition}
$$
\end{lemma}

\begin{proof}
Fix $\alpha,\beta\in F$ and put $\delta=\alpha\beta^{-1}$.
By~\eqref{eq:coefficient-partition},
$$
\begin{aligned}
 \left(\sum_{a\in A}R_{\sigma}\circ W_a\right)_{\alpha\beta}
 &=r_{\sigma}(\delta)\sum_{a\in A}\omega(\delta)^a\\
 &=r_{\sigma}(\delta)P_A(\omega(\delta))\\
 &=\1_{\{1\}}(\delta)=\1_{\{\alpha=\beta\}}.
\end{aligned}
$$
Hence~\eqref{eq:matrix-partition} holds.
\end{proof}

\begin{lemma}
\label{lem:matrix-balance}
For every finite nonempty $F\subset\Gamma$,
$$
 R\circ W_M+\sum_{a\in A\setminus\{0,M\}}K_a\circ R\circ W_a
   =R_{\sigma}\circ W_M+\sum_{a\in A\setminus\{0,M\}}\frac aM R_{\sigma}\circ W_a.
\numberthis{eq:matrix-balance}
$$
\end{lemma}

\begin{proof}
Fix $\alpha,\beta\in F$ and put $\delta=\alpha\beta^{-1}$.
By~\eqref{eq:coefficient-balance},
$$
\begin{aligned}
 &\left(R\circ W_M+
   \sum_{a\in A\setminus\{0,M\}}K_a\circ R\circ W_a\right)_{\alpha\beta}\\
 &\quad=r(\delta)\left[\omega(\delta)^M+
   \sum_{a\in A\setminus\{0,M\}}
   k_{a,M}(\ell(\delta))\omega(\delta)^a\right]\\
 &\quad=r_{\sigma}(\delta)\left[\omega(\delta)^M+
   \sum_{a\in A\setminus\{0,M\}}\frac aM\omega(\delta)^a\right]\\
 &\quad=\left(R_{\sigma}\circ W_M+
   \sum_{a\in A\setminus\{0,M\}}\frac aM R_{\sigma}\circ W_a
   \right)_{\alpha\beta}.
\end{aligned}
$$
Since $\alpha,\beta$ are arbitrary, this proves~\eqref{eq:matrix-balance}.
\end{proof}

\begin{lemma}
\label{lem:master-identity}
For every finite nonempty $F\subset\Gamma$,
$$
 R\circ W_M
 +\sum_{a\in A\setminus\{0,M\}}K_a\circ R\circ W_a+R_{\sigma}
 +\sum_{a\in A\setminus\{0,M\}}\left(1-\frac aM\right)R_{\sigma}\circ W_a
   =I_F.
\numberthis{eq:master-identity}
$$
\end{lemma}

\begin{proof}
Using~\eqref{eq:matrix-partition}, $R_{\sigma}\circ W_0=R_{\sigma}$, and
\eqref{eq:matrix-balance}, we obtain
$$
\begin{aligned}
 I_F
 &=R_{\sigma}+R_{\sigma}\circ W_M+
   \sum_{a\in A\setminus\{0,M\}}R_{\sigma}\circ W_a\\
 &=R\circ W_M+
   \sum_{a\in A\setminus\{0,M\}}K_a\circ R\circ W_a+R_{\sigma}\\
 &\quad+\sum_{a\in A\setminus\{0,M\}}
   \left(1-\frac aM\right)R_{\sigma}\circ W_a,
\end{aligned}
$$
which is~\eqref{eq:master-identity}.
\end{proof}

\subsection{A matrix deficit estimate and anchor preservation}
\label{sec:filter-deficit}

For Hermitian matrices, $B\preceq C$ means $C-B\succeq0$.  The positive terms in~\eqref{eq:master-identity} yield a comparison of the original and filtered matrices.

\begin{proposition}
\label{prop:matrix-deficit}
For every finite nonempty $F\subset\Gamma$,
$$
 R\preceq(1-c_A)I_F+c_AR_{\sigma},
 \qquad\text{equivalently}\qquad
 I_F-R\succeq c_A(I_F-R_{\sigma}).
\numberthis{eq:matrix-deficit}
$$
\end{proposition}

\begin{proof}
Recall that $R,R_\sigma,W_a\succeq0$, as established
after~\eqref{eq:matrix-WK}, and $K_a\succeq0$
by~\cref{lem:hyperbolic}.  Since $1-a/M\ge0$,
\cref{thm:schur-product} shows that every summand on the left-hand
side of~\eqref{eq:master-identity} is positive semidefinite.
Using $1-a/M\ge c_A$ and~\eqref{eq:matrix-partition}, we obtain
$$
\begin{aligned}
 I_F-R\circ W_M
 &\succeq R_\sigma+
   c_A\sum_{a\in A\setminus\{0,M\}}R_\sigma\circ W_a\\
 &=c_A(I_F-R_\sigma\circ W_M)+(1-c_A)R_\sigma\\
 &\succeq c_A(I_F-R_\sigma\circ W_M).
\end{aligned}
$$
Let $D=\diag(\omega(\alpha)^M)_{\alpha\in F}$.  Then
$$
 D^*D=I_F,\qquad
 R\circ W_M=DRD^*,\qquad R_\sigma\circ W_M=DR_\sigma D^*.
$$
Conjugating the last inequality by $D^*$ gives
$I_F-R\succeq c_A(I_F-R_\sigma)$, which is~\eqref{eq:matrix-deficit}.
\end{proof}

The next lemma transfers this finite matrix comparison to every point of $G$.

\begin{lemma}
\label{lem:filter}
Under the hypotheses~\eqref{eq:filter-assumptions} and~\eqref{eq:root-hypothesis},
the following hold:
\begin{enumerate}[label=(\roman*),leftmargin=2.2em]
\item For every $z\in G$,
$$
 0\le c_A\bigl(1-f_{\sigma}(z)\bigr)\le1-f(z).
\numberthis{eq:pointwise-deficit}
$$
\item The functions $f$ and $f_\sigma$ attain the value $1$ at exactly
the same points:
$$
 f(z)=1\quad\Longleftrightarrow\quad f_\sigma(z)=1
 \qquad\text{for all }z\in G.
$$
\item The sequence $y(n):=f_\sigma(g^n)$ belongs to $X_A$.
If $f(e_G)=1$, then $y\in X_A^{(1)}$.
\end{enumerate}
\end{lemma}

\begin{proof}
\textup{(i)} Choose a F\o lner sequence $(F_N)$ in $\Gamma$.  Explicitly, after
identifying $\Gamma$ with $(\Z/q\Z)\oplus\Z^\nu$, one may take the boxes
$(\Z/q\Z)\times\{-N,\ldots,N\}^\nu$.  Thus, for every
fixed $\delta\in\Gamma$,
$$
 \frac{|F_N\cap\delta^{-1}F_N|}{|F_N|}\longrightarrow1.
$$
This includes the finite-group case, where the whole group can be used
for every $N$.

For a fixed $z\in G$, set
$$
 p_{N,z}(\alpha)=\frac{\overline{z(\alpha)}}{\sqrt{|F_N|}}
 \quad(\alpha\in F_N).
$$
This is a unit vector.  Grouping the terms by
$\delta=\alpha\beta^{-1}$ gives
$$
\begin{aligned}
 p_{N,z}^*Rp_{N,z}
 &=\frac1{|F_N|}\sum_{\alpha,\beta\in F_N}
      r(\alpha\beta^{-1})z(\alpha\beta^{-1})\\
 &=\sum_{\delta\in\Gamma}r(\delta)z(\delta)
   \frac{|F_N\cap\delta^{-1}F_N|}{|F_N|}
   \longrightarrow f(z).
\end{aligned}
$$
The same calculation for $R_{\sigma}$ converges to $f_{\sigma}(z)$.
Both limits are justified by finite support.  Applying
\eqref{eq:matrix-deficit} to $p_{N,z}$ and passing to the limit yields
$$
 f(z)\le1-c_A+c_Af_{\sigma}(z).
$$
Rearranging and using $0\le f_{\sigma}\le1$ gives
\eqref{eq:pointwise-deficit}.  

\textup{(ii)} If $f(z_0)=1$, then
$0\ge c_A(1-f_{\sigma}(z_0))\ge0$, and $c_A>0$ gives $f_{\sigma}(z_0)=1$.
Conversely, if $f_{\sigma}(z_0)=1$, then~\eqref{eq:conditional-expectation}
gives
$$
 \int_H\bigl(1-f(z_0h)\bigr)\,dm_H(h)=0.
$$
The integrand is continuous and nonnegative, and Haar measure has full
support on $H$.  Hence $f(z_0h)=1$ for every $h\in H$; taking
$h=e_G$ gives $f(z_0)=1$.

\textup{(iii)} Since $0\le f_\sigma\le1$, we have $0\le y\le1$.
Evaluating the second identity of~\eqref{eq:filtered-partitions}
at $z=g^n$ gives
$$
 \sum_{a\in A}y(n-a)=\sum_{a\in A}f_\sigma(g^{n-a})=1
 \qquad\text{for all }n\in\Z.
$$
Thus $y\in X_A$.  If $f(e_G)=1$, part~\textup{(ii)} gives
$y(0)=f_\sigma(e_G)=1$, so $y\in X_A^{(1)}$.
\end{proof}

\subsection{An alternative proof by harmonic extension}
\label{sec:filter-harmonic}

In this section, we prove~\cref{lem:filter} by applying the strong minimum principle
to a harmonic extension of $f$ on $\R\times(0,M)$.

For a finite exponential sum
$F(\tau)=\sum_{q\in Q}d_qe^{iq\tau}$, where $Q\subset\R$ is finite,
its \emph{mean value} (or \emph{Bohr mean}) is
$$
 \mathcal M(F):=\lim_{T\to\infty}\frac1{2T}
                    \int_{-T}^T F(\tau)\,d\tau.
$$
Termwise integration shows that this limit exists and
$\mathcal M(F)=d_0$, with $d_0=0$ if $0\notin Q$.

\begin{proof}[Alternative proof of \cref{lem:filter}]
Let $h_\tau\in G$ be given by
$h_\tau(\gamma)=e^{i\tau\ell(\gamma)}$.  Since
$\mathcal M(e^{it\tau})=\1_{\{0\}}(t)$, finite support gives
$$
 f_\sigma(z)=\mathcal M\bigl(\tau\mapsto f(zh_\tau)\bigr)
 \in[0,1]\qquad(z\in G).
\numberthis{eq:harmonic-filter-average}
$$

For $0\le s\le M$, let
$$
 k_{s,M}(t)=
 \begin{cases}
  \dfrac{\sinh(st)}{\sinh(Mt)},&t\ne0,\\[6pt]
  \dfrac{s}{M},&t=0.
 \end{cases}
$$
In particular,
$$
 k_{0,M}(t)=0,\qquad k_{M,M}(t)=1\qquad(t\in\R).
$$
Define $V:G\times[0,M]\to\C$ by
$$
 V(w,s):=\sum_{\gamma\in\Gamma}r(\gamma)w(\gamma)
                           k_{s,M}(\ell(\gamma)).
\numberthis{eq:harmonic-auxiliary}
$$
Real-valuedness of $f$ implies
$r(\gamma^{-1})=\overline{r(\gamma)}$.  Since
$\ell(\gamma^{-1})=-\ell(\gamma)$ and $k_{s,M}$ is real and even,
$V$ is real-valued and continuous.  Moreover,
$$
 \partial_s^2 k_{s,M}(t)=t^2k_{s,M}(t),\qquad
 k_{s,M}(0)=\frac{s}{M},\qquad V(w,0)=0,\qquad V(w,M)=f(w).
$$
Consequently, for each $z\in G$, the function
$U_z(\tau,s):=V(zh_\tau,s)$ has the following properties:
\begin{enumerate}[label=(\roman*),leftmargin=2.2em]
\item $U_z$ is harmonic on $\R\times(0,M)$:
$$
 (\partial_\tau^2+\partial_s^2)U_z(\tau,s)=0
 \qquad\bigl((\tau,s)\in\R\times(0,M)\bigr).
$$
\item $U_z(\tau,0)=0$ for every $\tau\in\R$.
\item $U_z(\tau,M)=f(zh_\tau)\ge0$ for every $\tau\in\R$.
\end{enumerate}

We first show that $V\ge0$.  If
$$
 \mu:=\min_{(w,s)\in G\times[0,M]}V(w,s)
      =V(w_0,s_0)<0,
$$
then $0<s_0<M$ by the boundary conditions.  Since
$U_{w_0}\ge\mu$ and $U_{w_0}(0,s_0)=\mu$, the strong minimum
principle gives $U_{w_0}\equiv\mu$ on $\R\times(0,M)$.
Continuity at $s=0$ then gives $\mu=0$, a contradiction.

The coefficient identity~\eqref{eq:coefficient-balance}, with
$k_{a,M}(0)=a/M$, yields
$$
 \begin{aligned}
 \sum_{a\in A}V(z\omega^a,a)
 &=\sum_{\gamma\in\Gamma}r(\gamma)z(\gamma)
       \sum_{a\in A}k_{a,M}(\ell(\gamma))\omega(\gamma)^a\\
 &=\sum_{\gamma\in\Gamma}r_\sigma(\gamma)z(\gamma)
       \sum_{a\in A}\frac aM\omega(\gamma)^a
  =\sum_{a\in A}\frac aM f_\sigma(z\omega^a).
 \end{aligned}
\numberthis{eq:harmonic-function-balance}
$$
Combining this with the first partition
in~\eqref{eq:filtered-partitions} and $V(w,M)=f(w)$, we obtain
$$
\begin{aligned}
 1-f(z\omega^M)
 &=\sum_{a\in A\setminus\{M\}}V(z\omega^a,a)
   +\sum_{a\in A\setminus\{M\}}
       \left(1-\frac aM\right)f_\sigma(z\omega^a)\\
 &\ge c_A\sum_{a\in A\setminus\{M\}}f_\sigma(z\omega^a)\\
 &=c_A\bigl(1-f_\sigma(z\omega^M)\bigr)\ge0.
\end{aligned}
$$
Here $V\ge0$, $f_\sigma\ge0$, and $1-a/M\ge c_A$ for $a<M$.
Replacing $z$ by $z\omega^{-M}$ proves part~\textup{(i)}.

Part~\textup{(i)} implies $f(z)=1\Rightarrow f_\sigma(z)=1$.
Conversely, if $f_\sigma(z)=1$, let $D(\tau)=1-f(zh_\tau)$.
Then $0\le|D|^2\le D$ and
$\mathcal M(D)=1-f_\sigma(z)=0$ by~\eqref{eq:harmonic-filter-average}.
Writing $D(\tau)=\sum_{q\in Q}d_qe^{iq\tau}$ with $Q\subset\R$
finite and its frequencies distinct gives
$$
 0\le\sum_{q\in Q}|d_q|^2
   =\mathcal M(|D|^2)\le\mathcal M(D)=0.
$$
Hence $D\equiv0$ and $f(z)=1$, proving part~\textup{(ii)}.
For $y(n)=f_\sigma(g^n)$, the second
partition in~\eqref{eq:filtered-partitions} gives
$$
 0\le y\le1,\qquad
 \sum_{a\in A}y(n-a)=1\qquad(n\in\Z).
$$
If $f(e_G)=1$, part~\textup{(ii)} gives $y(0)=1$.
This proves part~\textup{(iii)}.
\end{proof}

Together with the Galois iteration in~\cref{sec:main-proof}, this
gives an alternative proof of the main theorems.

\subsection{The one-step filter}

Applying \cref{lem:filter} to $r=r_x$ gives the following conclusion.

\begin{proposition}\label{prop:one-filter}
Let $A\subset\Z$ be finite, nonempty, and normalized as in~\eqref{eq:normalization}, let $x\in X_A^{(1)}$ have the expansion~\eqref{eq:expansion-intro}, and let $K$ be the splitting field of $P_A$ over $\Q$.  For every $\sigma\in\Gal(K/\Q)$, the filtered sequence
$$
 y_\sigma(n):=\frac1{|A|}+\sum_{\substack{\lambda\in\Lambda_x\\
                                |\sigma(\lambda)|=1}}
                                c_\lambda\lambda^n
$$
belongs to $X_A^{(1)}$.
\end{proposition}

\begin{proof}
The singleton case is immediate.  Otherwise, the filtered coefficients
form the function $r_{\sigma}$ of~\eqref{eq:filtered-coefficients}, and
$y_\sigma(n)=f_{\sigma}(g^n)$ by~\eqref{eq:filtered-coefficients} and~\eqref{eq:orbit-values}.  \cref{lem:filter} therefore gives $y_\sigma\in X_A^{(1)}$.
\end{proof}

\section{Proof of the main theorems}
\label{sec:main-proof}

\begin{proof}[Proof of \cref{thm:expansion}]
If $|A|=1$, then $x\equiv1$ and the conclusions are immediate.
Assume $|A|\ge2$ and normalize $A$ as in~\eqref{eq:normalization}.
The existence and uniqueness of~\eqref{eq:expansion-intro} for every
$x\in X_A$ follow from~\cref{lem:recurrence}.  It remains to prove the
assertions about the projection for arbitrary $x\in X_A$.  Construct $r_x$,
$\Gamma$, $G$, $g$, and $f_x$ as in~\cref{lem:compact-dual}.  Let $K$
be the splitting field of $P_A$ over $\Q$, and write
$$
 \Gal(K/\Q)=\{\sigma_1,\ldots,\sigma_L\}.
$$
Starting from $r^{(0)}=r_x$, define recursively
$$
 r^{(j)}(\gamma)=r^{(j-1)}(\gamma)
 \1_{\{|\sigma_j(\gamma)|=1\}},
 \qquad
 f^{(j)}(z)=\sum_{\gamma\in\Gamma}r^{(j)}(\gamma)z(\gamma).
\numberthis{eq:iterated-filter}
$$

We apply \cref{lem:filter} at each step, using either of its two proofs.
The averaging formula preserves $0\le f^{(j)}\le1$, and the support
only decreases, so the root conditions~\eqref{eq:root-hypothesis}
remain valid.  Hence, for every $j$,
$$
 0\le f^{(j)}\le1,
 \qquad n\mapsto f^{(j)}(g^n)\in X_A.
\numberthis{eq:iteration-properties}
$$
Moreover, part~\textup{(ii)} of~\cref{lem:filter} gives
$$
 f^{(j)}(z)=1\quad\Longleftrightarrow\quad f^{(j-1)}(z)=1
 \qquad(z\in G).
\numberthis{eq:iteration-maximum-points}
$$

Let $r^{(L)}$ be the final coefficient function.  If
$\gamma\in\supp r^{(L)}\setminus\{1\}$, then
$$
 |\sigma(\gamma)|=1\qquad\text{for every }\sigma\in\Gal(K/\Q).
\numberthis{eq:all-conjugates-circle}
$$
The number $\gamma$ is an algebraic integer because it is a root of the monic polynomial $P_A$.  Since $K/\Q$ is normal, every algebraic conjugate of $\gamma$ equals $\sigma(\gamma)$ for some $\sigma\in\Gal(K/\Q)$.  By~\eqref{eq:all-conjugates-circle} and Kronecker's classical theorem, recalled in~\cite[Theorem~2.11]{FuSong}, $\gamma$ is a root of unity.

Conversely, a root of unity has all Galois conjugates on the unit circle and is never removed by~\eqref{eq:iterated-filter}.  Thus $r^{(L)}$ is obtained from $r_x$ by retaining exactly the root-of-unity frequencies.  Define
$$
 y(n)=f^{(L)}(g^n).
$$
By~\eqref{eq:iteration-properties}, $y\in X_A$.  Its nonconstant
frequencies are finitely many roots of unity, so their orders have
a finite least common multiple $N$ (take $N=1$ if there are none), and
$$
 y(n+N)=y(n)\qquad(n\in\Z).
$$
Iterating~\eqref{eq:iteration-maximum-points} at $z=g^n$ yields
$$
 x(n)=1\quad\Longleftrightarrow\quad y(n)=1\qquad(n\in\Z).
$$
In particular, $x(0)=1$ implies $y(0)=1$.

This proves~\cref{thm:expansion}, and hence~\cref{thm:main}.
\end{proof}

\section{Consequences}
\label{sec:consequences}

We first record a consequence of the recurrence that does not require
the condition $x(0)=1$ or the Galois filter.

\begin{corollary}\label{cor:finite-values}
If $x\in X_A$ takes only finitely many values, then $x$ is periodic.
In particular, every $\{0,1\}$-valued weak tiling is periodic.
\end{corollary}

\begin{proof}
The singleton case is immediate.  Otherwise, put $u(n)=x(n)-1/|A|$
and $U_n=(u(n),\ldots,u(n-M+1))^{\mathsf T}$.
By~\cref{lem:recurrence}, $U_n=C_A^nU_0$ with $C_A$ invertible.
The vectors $U_n$ take only finitely many values, so
$U_{n_1}=U_{n_2}$ for some $n_1<n_2$.  With $p=n_2-n_1$,
invertibility gives $U_{n+p}=U_n$ for every $n\in\Z$, hence
$x(n+p)=x(n)$.  Uniqueness of~\eqref{eq:expansion-intro} then gives
$\lambda^p=1$ for every $\lambda\in\Lambda_x$.
\end{proof}

For arbitrary normalized weights, the original sequence can be
nonperiodic, but the period of the projected weight depends only on
the roots of the mask polynomial.

\begin{corollary}\label{cor:period-bound}
Let
$$
 N_A=\lcm\{\operatorname{ord}(\zeta):P_A(\zeta)=0,
          \ \zeta\text{ is a root of unity}\},
$$
where $\operatorname{ord}(\zeta)$ is the multiplicative order of $\zeta$,
with the convention $N_A=1$ if the set is empty.  Whenever
$X_A^{(1)}\ne\varnothing$, it contains an $N_A$-periodic point.
\end{corollary}

\begin{proof}
Every nonconstant frequency of the projected point in~\cref{thm:expansion} is a root of $P_A$ and a root of unity.  Its order therefore divides $N_A$.
\end{proof}

The original weight can be nonperiodic, but the positions where it
attains its maximum value $1$ retain this periodicity.

\begin{corollary}\label{cor:maximum-set}
Let $x\in X_A^{(1)}$, and let $y$ be its periodic projection
in~\eqref{eq:ru-projection-intro}.  Then
$$
 \{n\in\Z:x(n)=1\}=\{n\in\Z:y(n)=1\}.
$$
In particular, $\{n\in\Z:x(n)=1\}$ is a union of residue classes
modulo $N_A$.
\end{corollary}

\begin{proof}
The assertion is immediate if $|A|=1$.  Otherwise, consider the
functions $f^{(j)}$ in~\eqref{eq:iterated-filter}.  By~\cref{lem:filter},
each filtering step preserves the set of points where the function
equals $1$.  Thus, for every $z\in G$,
$$
 f^{(0)}(z)=1\quad\Longleftrightarrow\quad f^{(L)}(z)=1.
$$
Since $x(n)=f^{(0)}(g^n)$ and $y(n)=f^{(L)}(g^n)$, the two sets
in the statement coincide.  The frequencies of $y$ have orders
dividing $N_A$, so $y$ is $N_A$-periodic, proving the final assertion.
\end{proof}

\begin{corollary}\label{cor:cyclotomic-obstruction}
If $|A|>1$ and none of the roots of $P_A$ is a root of unity, then $X_A^{(1)}$ is empty.
\end{corollary}

\begin{proof}
If $x\in X_A^{(1)}$, \cref{thm:main} produces a periodic $y\in X_A^{(1)}$.
By~\cref{lem:recurrence}, every nonconstant frequency of $y$ is a root
of $P_A$; by periodicity and uniqueness of the expansion, it is also
a root of unity.  Under the hypothesis there are no such frequencies,
so $y\equiv1/|A|$, contradicting $y(0)=1$.
\end{proof}

\begin{corollary}\label{cor:irreducible-mask}
Suppose $|A|>1$ and $P_A$ is irreducible over $\Q$.  If $X_A^{(1)}$ is nonempty, then $P_A$ is cyclotomic and every $x\in X_A$ is periodic.
\end{corollary}

\begin{proof}
By \cref{cor:cyclotomic-obstruction}, $P_A$ has a root of unity as a zero.  Its minimal polynomial is cyclotomic and divides $P_A$; since $P_A$ is monic and irreducible, the two polynomials coincide.
By~\cref{lem:recurrence}, all nonconstant frequencies of every $x\in X_A$ are roots of $P_A$.  They are therefore roots of unity, so every such $x$ is periodic.
\end{proof}

The next observation is another consequence of the recurrence alone.  It gives a simple range in which even unnormalized weights are periodic.

\begin{corollary}\label{cor:few-terms}
If $1\le|A|\le4$, then every $\lambda\in\C$ satisfying $P_A(\lambda)=0$ and $|\lambda|=1$ is a root of unity.  Consequently every $x\in X_A$ is periodic.
\end{corollary}

\begin{proof}
The case $|A|=1$ is immediate.  For $|A|=2$, after translation $P_A(z)=1+z^a$ with $a>0$, so every zero satisfies $z^{2a}=1$.
For $|A|=3$, write $P_A(z)=1+z^a+z^b$ with $0<a<b$.  If $P_A(\lambda)=0$ and $|\lambda|=1$, then $|1+\lambda^a|=1$, whence
$$
 \lambda^a+\lambda^{-a}=-1.
$$
Thus $\lambda^a$ is a nontrivial cube root of unity, and $\lambda^{3a}=1$.

For $|A|=4$, let $\lambda\in\C$ satisfy $P_A(\lambda)=0$ and $|\lambda|=1$.
Write the four summands of $P_A(\lambda)$ as $u_1,u_2,u_3,u_4$, so $|u_j|=1$ and $\sum_ju_j=0$.
If $u_1+u_2=0$, an opposite pair is already present.  Otherwise the unordered pair of points on the unit circle with sum $u_1+u_2$ is uniquely determined by its nonzero midpoint.  Since
$$
 u_1+u_2=(-u_3)+(-u_4),
$$
we obtain $\{u_1,u_2\}=\{-u_3,-u_4\}$.  In either case $\lambda^a=-\lambda^b$ for distinct $a,b\in A$, and hence $\lambda^{2|a-b|}=1$.
The periodicity assertion now follows from \cref{lem:recurrence}.
\end{proof}

\section{Applications}
\label{sec:applications}

The finite-valued consequence recovers ordinary tiling periodicity.
The projection theorem has a stronger conclusion for positive-definite
weights, which in turn yields rationality of finite spectra.

\subsection{Translational tilings}

A finite nonempty set $A\subset\Z$ tiles $\Z$ by translations if
there is a set $T\subset\Z$ such that every integer has exactly one
representation $n=a+t$ with $a\in A$ and $t\in T$.
Equivalently, $\1_A*\1_T\equiv1$.  The indicator $\1_T$ belongs to
$X_A$ and takes only two values, so~\cref{cor:finite-values} makes
$T$ periodic.  This is Newman's classical theorem~\cite{Newman}.
Here periodicity concerns the given tiling complement, whereas
\cref{thm:main} supplies a periodic representative of a possibly
nonperiodic normalized weak tiling.

\subsection{Positive-definite weights}

A sequence $x:\Z\to\C$ is \emph{positive definite} if
$\bigl(x(j-k)\bigr)_{j,k\in E}\succeq0$ for every finite
$E\subset\Z$.  This is an additional condition on $x$; positive
definiteness of the coefficient function $r_x$ on $\Gamma$
in~\cref{lem:compact-dual} follows already from $x\ge0$.
The finite-group version of a normalized positive-definite weak tiling
is used in~\cite{KMMS,KLMS}.

\begin{corollary}\label{cor:positive-definite}
Let $x\in X_A^{(1)}$.  If all coefficients $c_\lambda$ in
\eqref{eq:expansion-intro} are nonnegative real numbers, then $x$ is
$N_A$-periodic.  In particular, this holds whenever $x$ is positive
definite on $\Z$.
\end{corollary}

\begin{proof}
Let $y$ be the projection in~\eqref{eq:ru-projection-intro}.
By~\cref{thm:expansion},
$$
 0=x(0)-y(0)=
 \sum_{\substack{\lambda\in\Lambda_x\\
                 \lambda\text{ is not a root of unity}}}c_\lambda.
$$
If the coefficients are nonnegative, this sum can vanish only when
there are no frequencies that are not roots of unity.  Thus $x=y$,
and every remaining frequency has order dividing $N_A$.

If $x$ is positive definite, its coefficients are nonnegative.
Indeed, for $\lambda\in\Lambda_x$ and $N\ge1$, testing the matrix
$\bigl(x(j-k)\bigr)_{0\le j,k<N}$ against
$N^{-1}(1,\lambda,\ldots,\lambda^{N-1})^{\mathsf T}$ gives
$$
\begin{aligned}
 0\le\frac1{N^2}\sum_{j,k=0}^{N-1}
       \lambda^{k-j}x(j-k)
 &=\sum_{\gamma\in\{1\}\cup\Lambda_x}r_x(\gamma)
   \left|\frac1N\sum_{j=0}^{N-1}(\gamma/\lambda)^j\right|^2\\
 &\longrightarrow r_x(\lambda)=c_\lambda.
\end{aligned}
$$
The limit follows from the finite support and the elementary averages
of distinct unit-modulus exponentials.
\end{proof}

\subsection{Rationality of finite spectra}

Galois methods for spectra of finite integer sets already appear in
Konyagin and \L aba~\cite[Section~3]{KonyaginLaba}, where orthogonality
and Galois averaging are combined under an irreducibility assumption.
The general rationality question is recorded
in~\cite[Question~5]{KolountzakisSurvey}.

Let $A\subset\Z$ be finite and nonempty, and put $d=|A|$.
A finite set $\Lambda\subset\R$ is a \emph{spectrum} of $A$ if the
functions $a\mapsto e^{2\pi i\lambda a}$, $\lambda\in\Lambda$, form
an orthogonal basis of $L^2(A)$ with counting measure.
In this case $|\Lambda|=d$, the elements of $\Lambda$ are distinct
modulo $1$, and
$$
 \sum_{a\in A}e^{2\pi i(\lambda-\mu)a}
 =d\,\1_{\{\lambda=\mu\}}\qquad(\lambda,\mu\in\Lambda).
$$
Replacing $\Lambda$ by a translate preserves spectrality, so we may
normalize it by $0\in\Lambda$.

Recall from~\cref{cor:period-bound} that
$$
 N_A=\lcm\{\operatorname{ord}(\zeta):P_A(\zeta)=0,
          \ \zeta\text{ is a root of unity}\},
$$
where $\operatorname{ord}(\zeta)$ is the multiplicative order of $\zeta$,
with $N_A=1$ if the set is empty.

\begin{corollary}\label{cor:rational-spectrum}
If $\Lambda$ is a spectrum of $A$ and $0\in\Lambda$, then
$\Lambda\subset N_A^{-1}\Z$.  In particular, $\Lambda\subset\Q$.
\end{corollary}

\begin{proof}
Define the spectral weight
$$
 w(n)=\frac1{d^2}
       \left|\sum_{\lambda\in\Lambda}e^{2\pi i\lambda n}\right|^2
     =\frac1{d^2}\sum_{\lambda,\mu\in\Lambda}
       e^{2\pi i(\lambda-\mu)n}.
$$
Then $0\le w\le1$, $w(0)=1$, and orthogonality gives
$$
\begin{aligned}
 \sum_{a\in A}w(n-a)
 &=\frac1{d^2}\sum_{\lambda,\mu\in\Lambda}
   e^{2\pi i(\lambda-\mu)n}
   \sum_{a\in A}e^{-2\pi i(\lambda-\mu)a}=1.
\end{aligned}
$$
Thus $w\in X_A^{(1)}$.  After equal frequencies are grouped, every
coefficient is the number of corresponding ordered pairs divided by
$d^2$, hence is positive.  By~\cref{cor:positive-definite}, $w$ is
$N_A$-periodic.  Uniqueness of its finite exponential expansion yields
$$
 e^{2\pi iN_A(\lambda-\mu)}=1
 \qquad(\lambda,\mu\in\Lambda).
$$
Taking $\mu=0$ gives $N_A\lambda\in\Z$ for every $\lambda\in\Lambda$.
\end{proof}

The rationality conclusion recovers the finite-spectrum result of
Fu and Song~\cite[Theorem~1.2]{FuSong}.  The argument here isolates
why it follows from the more general projection theorem: nonnegative
frequency coefficients prevent cancellation among deleted terms.

\section{A nonperiodic normalized weak tiling}
\label{sec:example}

The positive-definiteness assumption in~\cref{cor:positive-definite}
cannot be omitted: a normalized weak tiling can be nonperiodic even
when its frequencies include roots of unity.  We give an example with
$\gcd(A-A)=1$, so the tile is not confined to a residue class modulo
any integer greater than $1$.

\begin{proposition}\label{prop:nonperiodic-example}
Let
$$
 B=\{0,6,7,8,14\},\qquad
 A=4B+\{0,1\}=\{0,1,24,25,28,29,32,33,56,57\}.
$$
Then no translate of $A$ is contained in a proper subgroup of $\Z$, and $X_A^{(1)}$ contains both a nonperiodic point and a $20$-periodic point.
\end{proposition}

\begin{proof}
We first find $\alpha\in\C$ such that $P_B(\alpha)=0$, $|\alpha|=1$, and $\alpha$ is not a root of unity.  With $t=z+z^{-1}$, direct expansion gives
$$
\begin{aligned}
 z^{-7}P_B(z)
 &=z^{-7}+z^{-1}+1+z+z^7\\
 &=t^7-7t^5+14t^3-6t+1\\
 &=(t^2+t-1)h(t),
\end{aligned}
$$
where
$$
 h(t)=t^5-t^4-5t^3+4t^2+5t-1.
$$
Since $h(0)=-1$ and $h(1)=3$, choose a zero $t_0\in(0,1)$ and set
$$
 \theta=\arccos(t_0/2),\qquad \alpha=e^{i\theta}.
$$
Then $\alpha+\alpha^{-1}=t_0$, so $P_B(\alpha)=P_B(\alpha^{-1})=0$.

To see that $\alpha$ is not a root of unity, reduce $h$ modulo $2$:
$$
 \overline h(t)=t^5+t^4+t^3+t+1\in\mathbb F_2[t].
$$
This polynomial has no root in $\mathbb F_2$, and its remainder modulo the only irreducible quadratic $t^2+t+1$ is $t+1\ne0$.  A reducible polynomial of degree $5$ has a factor of degree at most $2$, so $\overline h$ is irreducible.  Hence $h$ is irreducible over $\Q$ and $[\Q(t_0):\Q]=5$.
If $\alpha$ had order $N$, then $N>2$ and
$$
 [\Q(\alpha+\alpha^{-1}):\Q]=\varphi(N)/2,
$$
because $\Q(\alpha)$ is a quadratic extension of its real subfield $\Q(\alpha+\alpha^{-1})$.  Thus $\varphi(N)=10$, which forces $N=11$ or $22$.
The corresponding minimal polynomials of $\alpha+\alpha^{-1}$ are
$$
 t^5+t^4-4t^3-3t^2+3t+1,
 \qquad
 t^5-t^4-4t^3+3t^2+3t-1,
$$
respectively.  The first follows by expanding $1+\sum_{j=1}^5(\alpha^j+\alpha^{-j})=0$ for $N=11$.  If the first polynomial is denoted by $p_{11}(t)$, the second is $-p_{11}(-t)$.  Neither equals $h$, a contradiction.

Define
$$
 u(n)=\1_{5\Z}(n),\qquad
 v(n)=\frac15+\frac1{10}\cos(n\theta).
$$
The set $B$ is a complete residue system modulo $5$, so $u\in X_B$.  Also $1/10\le v\le3/10$, and
$$
 \sum_{b\in B}v(n-b)
 =1+\frac1{10}\operatorname{Re}\bigl(\alpha^nP_B(\alpha^{-1})\bigr)=1,
$$
so $v\in X_B$ as well.  Interleave these two sequences by setting
$$
 x(4n)=u(n),\qquad x(4n+2)=v(n),\qquad
 x(4n+1)=x(4n+3)=0.
$$
Then $0\le x\le1$ and $x(0)=1$.  For $j\in\{0,1,2,3\}$,
$$
 \sum_{a\in A}x(4n+j-a)
 =\begin{cases}
   \sum_{b\in B}u(n-b)=1,&j=0,1,\\
   \sum_{b\in B}v(n-b)=1,&j=2,3.
  \end{cases}
$$
Thus $x\in X_A^{(1)}$.  Since $0,1\in A$, every translate of $A$ contains two consecutive integers and cannot lie in a proper subgroup of $\Z$.

Finally, $x(m)=1$ holds exactly when $m\in20\Z$.  Any nonzero period $p$ would therefore be $p=20q$ with $q\ne0$, and would imply $v(n+5q)=v(n)$ for all $n$.  Taking $n=0$ gives $\cos(5q\theta)=1$, or $\alpha^{5q}=1$, a contradiction.  Hence $x$ is nonperiodic.

Replacing $v$ with the constant sequence $1/5$ gives
$$
 y(4n)=u(n),\qquad y(4n+2)=\frac15,\qquad
 y(4n+1)=y(4n+3)=0.
$$
The same calculation shows that $y\in X_A^{(1)}$, and $y$ is $20$-periodic.
To identify the projection, use the residue-class expansion
$$
 x(m)-y(m)=\frac1{80}\sum_{j=0}^3 i^{j(m-2)}
 \left(e^{i(m-2)\theta/4}+e^{-i(m-2)\theta/4}\right).
$$
Every frequency in this expression has fourth power $\alpha$ or
$\alpha^{-1}$, so none is a root of unity.  Since $y$ is periodic,
uniqueness of~\eqref{eq:expansion-intro} shows that it is exactly
the projection in~\eqref{eq:ru-projection-intro}.
\end{proof}

The coexistence of the two types of frequencies is also visible in the factorization
$$
 P_B(z)=\Phi_5(z)S(z),\qquad
 S(z)=z^{10}-z^9+z^5-z+1,
$$
where $\Phi_5(z)=1+z+z^2+z^3+z^4$ and $S(z)=z^5h(z+z^{-1})$.
Consequently,
$$
 P_A(z)=(1+z)\Phi_5(z^4)S(z^4).
$$
The first two factors contribute roots of unity, while $S(z^4)$ has roots $\lambda$ with $|\lambda|=1$ that are not roots of unity, as used in the construction.

\section{Dynamics of the shift on the space of weak tilings}
\label{sec:dynamics}

Let $A\subset\Z$ be finite and nonempty, and write $d=|A|$.
The exponential expansion in~\cref{lem:recurrence} gives a single
system of coordinates for all of $X_A$.  In these coordinates the
shift is an isometry, and the frequencies with nonzero coefficients
determine the rotation on each orbit closure.  Throughout this section,
invariant measures are understood to have total mass one.

\subsection{Coefficient coordinates and topological entropy}

Let $\mathcal{Z}_A=\{\lambda\in\T:P_A(\lambda)=0\}$.  Allowing zero coefficients in
\eqref{eq:expansion-intro}, we write
$$
 x(n)=\frac1d+\sum_{\lambda\in\mathcal{Z}_A}c_\lambda(x)\lambda^n,
 \qquad c_{\bar\lambda}(x)=\overline{c_\lambda(x)}.
\numberthis{eq:dynamics-coefficients}
$$
The constant weak tiling $x_*\equiv1/d$ corresponds to the zero vector.

\begin{proposition}\label{prop:shift-isometry}
Let $c(x)=(c_\lambda(x))_{\lambda\in\mathcal{Z}_A}$ and
$\mathcal C_A=c(X_A)$.  The following hold.
\begin{enumerate}[label=(\roman*),leftmargin=2.2em]
\item The map $c:X_A\to\mathcal C_A$ is an affine homeomorphism, and
$\mathcal C_A$ is a compact convex set of real dimension $|\mathcal{Z}_A|$.
\item Let $D:\mathcal C_A\to\mathcal C_A$ be the diagonal isometry
defined by $(Dv)_\lambda=\lambda v_\lambda$ for $\lambda\in\mathcal{Z}_A$.
Then it gives the topological conjugacy
$$
 c\colon (X_A,\sigma)\longrightarrow(\mathcal C_A,D),
 \qquad c\circ\sigma=D\circ c.
$$
In particular, $\rho(x,y)=\|c(x)-c(y)\|_2$ is a compatible metric on
$X_A$ for which $\sigma$ is an isometry.
\item $h_{\mathrm{top}}(X_A,\sigma)=0$.
\end{enumerate}
\end{proposition}

\begin{proof}
\textup{(i)} If $\mathcal{Z}_A=\varnothing$, then $X_A=\{x_*\}$ and
$\mathcal C_A=\{0\}$.  Assume $\mathcal{Z}_A\ne\varnothing$, so $d\ge2$.
Put $I\coloneqq\{0,\ldots,|\mathcal{Z}_A|-1\}$ and fix an ordering $\prec$ of
$\mathcal{Z}_A$.  Write $c(x)$ as a column vector in this order and let
$V\coloneqq(\lambda^n)_{n\in I,\,\lambda\in\mathcal{Z}_A}$.  By
\eqref{eq:dynamics-coefficients}, the $n$th component of the product is
$$
 (Vc(x))_n=\sum_{\lambda\in\mathcal{Z}_A}\lambda^n c_\lambda(x)
 =x(n)-1/d,\qquad n\in I.
$$
The Vandermonde determinant is
$\det V=\prod_{\lambda\prec\mu}(\mu-\lambda)\ne0$.
Thus $c(x)=V^{-1}(x(n)-1/d)_{n\in I}$ is continuous and affine.
Equation~\eqref{eq:dynamics-coefficients} gives injectivity.  Since
$X_A$ is compact, $c$ is a homeomorphism onto its image $\mathcal C_A$,
which is compact and convex.

Let $E=\{v\in\C^{\mathcal{Z}_A}:v_{\bar\lambda}=\overline{v_\lambda}
\text{ for all }\lambda\in\mathcal{Z}_A\}$.  Each real root contributes
one real dimension and each nonreal conjugate pair contributes two, so
$\dim_{\R}E=|\mathcal{Z}_A|$.  For $v\in E$, put
$x_v(n)=1/d+\sum_{\lambda\in\mathcal{Z}_A}v_\lambda\lambda^n$.
Then $x_v$ is real-valued and
$$
 (\1_A*x_v)(n)
 =1+\sum_{\lambda\in\mathcal{Z}_A}v_\lambda\lambda^nP_A(\lambda^{-1})
 =1 .
$$
If $\|v\|_1<1/d$, then $|x_v(n)-1/d|\le\|v\|_1$ for all $n$,
so $0<x_v(n)<2/d\le1$.  Hence $x_v\in X_A$ and $c(x_v)=v$.
It follows that $\{v\in E:\|v\|_1<1/d\}\subseteq\mathcal C_A\subseteq E$,
and therefore $\dim_{\R}\mathcal C_A=|\mathcal{Z}_A|$.

\textup{(ii)} By uniqueness in~\eqref{eq:dynamics-coefficients},
$c_\lambda(\sigma x)=\lambda c_\lambda(x)$ for every
$\lambda\in\mathcal{Z}_A$.  Hence $c\circ\sigma=D\circ c$.  Since
$|\lambda|=1$, we have
$$
 \rho(\sigma x,\sigma y)=\|Dc(x)-Dc(y)\|_2
 =\|c(x)-c(y)\|_2=\rho(x,y).
$$
Compatibility of $\rho$ follows from \textup{(i)}.

\textup{(iii)} By \textup{(ii)}, $\sigma$ is an isometry of the compact
metric space $(X_A,\rho)$ and hence has zero topological entropy.
\end{proof}

\begin{remark}
We record four consequences of~\cref{prop:shift-isometry}.
\begin{enumerate}[label=(\roman*),leftmargin=2.2em]
\item The space $X_A$ is path connected.  For $x,y\in X_A$, convexity
gives the path $t\mapsto(1-t)x+ty$, $0\le t\le1$, from $x$ to $y$.
\item The shift is equicontinuous and distal.  Indeed,
$\rho(\sigma^n x,\sigma^n y)=\rho(x,y)$ for every $n\in\Z$, so all
iterates preserve distances and distinct points stay a positive distance apart.
\item The unique fixed point is $x_*\equiv1/d$.
\item If $X_A\ne\{x_*\}$, the shift has no dense orbit and is not
topologically transitive.  Since $\sigma$ fixes $x_*$ and preserves $\rho$, the function
$$
 R(x)\coloneqq\rho(x,x_*)^2
 =\sum_{\lambda\in\mathcal{Z}_A}|c_\lambda(x)|^2
$$
is continuous and shift-invariant.  Choose $y\ne x_*$ and $0<a<R(y)$.
Then $\{R<a\}$ and $\{R>a\}$ are disjoint nonempty invariant open
sets, which rules out topological transitivity.
\end{enumerate}
\end{remark}

\subsection{Orbit closures and invariant measures}

Under the conjugacy of~\cref{prop:shift-isometry}, the shift acts on
$\mathcal C_A$ by multiplication by $g_A$, where
$$
 g_A\coloneqq(\lambda)_{\lambda\in\mathcal{Z}_A},
 \qquad H_A\coloneqq\overline{\{g_A^n:n\in\Z\}}\subseteq\T^{\mathcal{Z}_A}.
$$
Since $g_A^n\mathcal C_A=\mathcal C_A$ for all $n\in\Z$ and
$\mathcal C_A$ is closed, coordinatewise multiplication restricts to
the continuous action
$$
 H_A\times\mathcal C_A\longrightarrow\mathcal C_A,\qquad
 (h,c)\longmapsto hc\coloneqq(h_\lambda c_\lambda)_{\lambda\in\mathcal{Z}_A}.
$$
For $c\in\mathcal C_A$, its stabilizer is
$$
 H_{A,c}\coloneqq\{h\in H_A:hc=c\}
 =\{h\in H_A:h_\lambda=1\text{ whenever }c_\lambda\ne0\}.
$$
Let $m_{H_A}$ denote the normalized Haar measure on $H_A$.

\begin{proposition}\label{prop:orbit-measures}
For $c\in\mathcal C_A$, put $Y_c\coloneqq\overline{\{D^nc:n\in\Z\}}$
and let $p_c:H_A\to Y_c$ be the orbit map $p_c(h)=hc$.
\begin{enumerate}[label=(\roman*),leftmargin=2.2em]
\item Let $T_c$ denote translation by $g_AH_{A,c}$ on $H_A/H_{A,c}$.
The map $p_c$ induces the topological conjugacy
$$
 \Phi_c\colon(H_A/H_{A,c},T_c)\longrightarrow(Y_c,D|_{Y_c}),
 \qquad hH_{A,c}\longmapsto hc.
$$
\item The system $(Y_c,D|_{Y_c})$ is minimal and uniquely ergodic.
Its unique invariant measure
$$
 \mu_c\coloneqq(p_c)_*m_{H_A}
 =\int_{H_A}\delta_{hc}\,dm_{H_A}(h)
$$
has discrete spectrum.
\item The ergodic $D$-invariant probability measures on $\mathcal C_A$
are precisely $\{\mu_c:c\in\mathcal C_A\}$.  Every $D$-invariant
probability measure $\mu$ on $\mathcal C_A$ satisfies
$$
 \mu=\int_{\mathcal C_A}\mu_c\,d\mu(c).
$$
\end{enumerate}
\end{proposition}

\begin{proof}
\textup{(i)} Since $D^nc=g_A^nc$ and $\langle g_A\rangle$ is dense
in $H_A$, compactness gives $Y_c=H_Ac$.  Moreover,
$hc=kc$ if and only if $k^{-1}h\in H_{A,c}$.  Thus $\Phi_c$ is a
continuous bijection from a compact space onto the Hausdorff space $Y_c$,
hence a homeomorphism.  For every $h\in H_A$,
$$
 \Phi_c(T_c(hH_{A,c}))=(g_Ah)c=D\Phi_c(hH_{A,c}).
$$
Projection onto the coordinates with $c_\lambda\ne0$ identifies
$H_A/H_{A,c}$ with the dual of $\langle\lambda:c_\lambda\ne0\rangle$,
as in~\cref{sec:preliminaries}.  In particular, the quotient is a
product of a finite cyclic group and a torus.

\textup{(ii)} The coset $g_AH_{A,c}$ generates a dense subgroup of
$H_A/H_{A,c}$, so $T_c$ is a minimal rotation.  By the structure theorem
for equicontinuous systems, this rotation is uniquely ergodic with
normalized Haar measure and has discrete spectrum~\cite[Theorem~2.11,
Proposition~2.12 and Section~2.3]{HermleKreidler}.  The quotient map
sends $m_{H_A}$ to normalized Haar measure on $H_A/H_{A,c}$, so
\textup{(i)} gives $\mu_c=(p_c)_*m_{H_A}$.

\textup{(iii)} Each $\mu_c$ is ergodic by \textup{(ii)}.  Conversely,
an ergodic measure $\nu$ has a generic point $c$, so $\nu(Y_c)=1$.
Uniqueness in \textup{(ii)} gives $\nu=\mu_c$.
For the last assertion, $D$-invariance implies $H_A$-invariance by
density of $\langle g_A\rangle$.  Hence for every $f\in C(\mathcal C_A)$,
$$
 \int_{\mathcal C_A}f\,d\mu
 =\int_{H_A}\int_{\mathcal C_A} f(hc)\,d\mu(c)\,dm_{H_A}(h)
 =\int_{\mathcal C_A}\left(\int_{\mathcal C_A}f\,d\mu_c\right)d\mu(c),
$$
by Fubini's theorem and the definition of $\mu_c$.
\end{proof}

\begin{corollary}\label{cor:ergodic-measure-count}
If $\mathcal{Z}_A=\varnothing$, then $\delta_{x_*}$ is the unique
invariant measure on $X_A$.  If $\mathcal{Z}_A\ne\varnothing$,
then $(X_A,\sigma)$ has exactly $2^{\aleph_0}$ ergodic measures.
\end{corollary}

\begin{proof}
If $\mathcal{Z}_A=\varnothing$, then~\eqref{eq:dynamics-coefficients}
gives $X_A=\{x_*\}$.  Otherwise, fix $\lambda\in\mathcal{Z}_A$ and set
$$
 x_t(n)\coloneqq\frac1d+t\operatorname{Re}(\lambda^n),
 \qquad 0<t<1/d.
\numberthis{eq:single-frequency-family}
$$
Since $P_A(\lambda^{-1})=0$ and $0<x_t(n)<2/d\le1$, we have $x_t\in X_A$.
By uniqueness of the coefficients, $c(x_t)=tv$ for a fixed $v\ne0$.
Thus $R(x_t)=t^2\|v\|_2^2$ is strictly increasing in $t$.  As $R$ is
continuous and shift-invariant, the orbit closures of the $x_t$ are
pairwise disjoint.  By~\cref{prop:orbit-measures}, each carries a distinct
ergodic measure.  Conversely, every ergodic measure is $\mu_c$ for some
$c\in\mathcal C_A$, and $\mathcal C_A\subseteq\C^{\mathcal{Z}_A}$ has
cardinality at most $2^{\aleph_0}$ since $\mathcal{Z}_A$ is finite.
\end{proof}

\subsection{Periodic points and the periodic projection}

Recall the integer $N_A$ from~\cref{cor:period-bound}, the least common
multiple of the orders of the roots of unity in $\mathcal{Z}_A$, with
$N_A=1$ if there are none.  By uniqueness in~\eqref{eq:dynamics-coefficients},
$\sigma^p x=x$ for $p\ge1$ precisely when $\lambda^p=1$ for each
$\lambda$ with $c_\lambda(x)\ne0$.  It follows that
$$
 \operatorname{Per}(\sigma)=\operatorname{Fix}(\sigma^{N_A})
 =\{x\in X_A:c_\lambda(x)=0\text{ whenever }\lambda
       \text{ is not a root of unity}\}.
$$
By~\cref{thm:expansion}, retaining just the root-of-unity coordinates
defines a map $\pi:X_A\to X_A$ given by
$$
 \pi(x)(n)=\frac1d+
 \sum_{\substack{\lambda\in\mathcal{Z}_A\\\lambda^{N_A}=1}}
       c_\lambda(x)\lambda^n .
$$
This coordinate projection is continuous and affine, commutes with the
shift, and fixes every periodic point.  Thus it is a retraction onto
$\operatorname{Per}(\sigma)$.

\begin{proposition}\label{prop:periodic-orbits-measures}
The following hold for $(X_A,\sigma)$.
\begin{enumerate}[label=(\roman*),leftmargin=2.2em]
\item If every $\lambda\in\mathcal{Z}_A$ is a root of unity, then
$\sigma^{N_A}=\mathrm{id}_{X_A}$.  The ergodic measures are precisely
the uniform measures on periodic orbits.
\item If some $\lambda\in\mathcal{Z}_A$ is not a root of unity, then
$\operatorname{Per}(\sigma)$ is a proper closed subset of $X_A$ and
$(X_A,\sigma)$ has $2^{\aleph_0}$ nonatomic ergodic measures.
\end{enumerate}
\end{proposition}

\begin{proof}
\textup{(i)} Since $\lambda^{N_A}=1$ for all $\lambda\in\mathcal{Z}_A$,
\eqref{eq:dynamics-coefficients} gives $\sigma^{N_A}x=x$ for every $x\in X_A$.
Each orbit is therefore finite and carries its uniform measure.
By~\cref{prop:orbit-measures} these are precisely the ergodic measures.

\textup{(ii)} Fix such a $\lambda$ and take $x_t$ as
in~\eqref{eq:single-frequency-family}.  Its coefficient at $\lambda$
is $t/2\ne0$, so a period $p\ge1$ would imply $\lambda^p=1$.
Thus $x_t\notin\operatorname{Per}(\sigma)=\operatorname{Fix}(\sigma^{N_A})$,
which proves that this closed set is proper.  As in the proof
of~\cref{cor:ergodic-measure-count}, the orbit closures of the $x_t$
are pairwise disjoint.  Each is conjugate to a rotation on an infinite
compact group by~\cref{prop:orbit-measures}.  Its Haar measure is
nonatomic, since all singletons have equal mass.  This gives
$2^{\aleph_0}$ distinct nonatomic ergodic measures, and the matching
upper bound follows from~\cref{cor:ergodic-measure-count}.
\end{proof}

\subsection{Weak tilings normalized at the origin}

Set $m_A=\max_{x\in X_A}x(0)$, which exists by compactness.  Shift
invariance gives
$$
 m_A=\sup_{x\in X_A,\,n\in\Z}x(n).
$$
Thus $X_A^{(1)}=\varnothing$ if and only if $m_A<1$, or equivalently,
if and only if $x(n)<1$ for every $x\in X_A$ and $n\in\Z$.
By~\cref{thm:main}, $X_A^{(1)}\ne\varnothing$ if and only if there
exists a periodic $x\in X_A$ with $x(0)=1$.  For $d>1$ such a point is nonconstant.  One of its
active frequencies is therefore a root of unity, and
\eqref{eq:single-frequency-family} for that frequency gives continuum
many ergodic periodic-orbit measures.

\begin{proposition}\label{prop:normalized-weak-tilings}
Let $Y_A\coloneqq\{x\in X_A:x(n)=1\text{ for some }n\in\Z\}$.
\begin{enumerate}[label=(\roman*),leftmargin=2.2em]
\item The set $X_A^{(1)}$ is closed and convex, and
$\sigma^{N_A}X_A^{(1)}=X_A^{(1)}$.
\item The set $Y_A$ is closed and shift-invariant, and
$$
 Y_A=\bigcup_{j=0}^{N_A-1}\sigma^{-j}X_A^{(1)}.
$$
\item For $x\in X_A$, let $\mu_x$ denote the unique invariant measure
on its orbit closure.  Then
$$
 \mu_x(X_A^{(1)})=
 \frac{1}{N_A}\#\{0\le j<N_A:x(j)=1\}.
\numberthis{eq:normalized-set-mass}
$$
\end{enumerate}
\end{proposition}

\begin{proof}
\textup{(i)} The evaluation $x\mapsto x(0)$ is continuous and affine,
so its level set $X_A^{(1)}$ is closed and convex.
By~\cref{thm:expansion}, $x(n)=1$ if and only if $\pi(x)(n)=1$.
Since $\pi(x)$ is $N_A$-periodic, $x(0)=1$ if and only if $x(N_A)=1$.
This gives $\sigma^{N_A}X_A^{(1)}=X_A^{(1)}$.

\textup{(ii)} Since $x(n)=1$ if and only if $\sigma^nx\in X_A^{(1)}$,
we have $Y_A=\bigcup_{n\in\Z}\sigma^{-n}X_A^{(1)}$.  This union is
shift-invariant.  By \textup{(i)} it reduces to the union over
$0\le n<N_A$ and hence is closed.

\textup{(iii)} The map $\pi$ carries the orbit closure of $x$ onto the finite orbit
of $\pi(x)$.  Its pushforward of $\mu_x$ is invariant, hence uniform
on that orbit.  The level-set equality at $n=0$ gives
$X_A^{(1)}=\pi^{-1}(X_A^{(1)}\cap\operatorname{Per}(\sigma))$.
Averaging the indicator of $X_A^{(1)}$ over the $N_A$ iterates of
$\pi(x)$ now gives~\eqref{eq:normalized-set-mass}.
\end{proof}

\begin{corollary}\label{cor:normalized-positive-mass}
The set $X_A^{(1)}$ is nonempty if and only if there exists an ergodic
measure $\mu$ on $X_A$ with $\mu(X_A^{(1)})>0$.  In this case $\mu$
can be chosen to be the uniform measure on a periodic orbit.
\end{corollary}

\begin{proof}
If $X_A^{(1)}\ne\varnothing$, then~\cref{thm:main} gives a periodic
point $y\in X_A^{(1)}$.  The uniform measure $\mu_y$ on its orbit is
ergodic, and $y(0)=1$ gives $\mu_y(X_A^{(1)})\ge 1/N_A>0$
by~\eqref{eq:normalized-set-mass}.  Conversely, $\mu(X_A^{(1)})>0$
implies $X_A^{(1)}\ne\varnothing$.
\end{proof}

\subsection{Normalization under topological conjugacy}
\label{sec:normalization-conjugacy}
Take $A_+=\{0,1\}$ and $A_-=\{0,1,2,5\}$.  Both masks have $-1$
as their only unit-circle zero.  For $A_-$, a unit-circle zero of
$P(z)=1+z+z^2+z^5$ must also be a zero of $P^*(z)=z^5P(z^{-1})$.
Since $P(z)-P^*(z)=z(1-z)(1+z)^2$ and $P(1)=4$, only $-1$ is
possible, and $P(-1)=0$.  Thus~\eqref{eq:dynamics-coefficients} gives
$$
\begin{aligned}
 X_{A_+}&=\{x_c(n)=1/2+c(-1)^n: |c|\le1/2,\ c\in\R\},\\
 X_{A_-}&=\{y_c(n)=1/4+c(-1)^n: |c|\le1/4,\ c\in\R\}.
\end{aligned}
$$
In each case the shift sends $c$ to $-c$, and $x\mapsto x/2$ is a
topological conjugacy from $X_{A_+}$ to $X_{A_-}$.  Yet $x_{1/2}(0)=1$,
whereas every coordinate in $X_{A_-}$ is at most $1/2$.  Thus
$X_{A_+}^{(1)}\ne\varnothing$ and $X_{A_-}^{(1)}=\varnothing$,
although the shifts are topologically conjugate.

The shift on $X_{A_-}$ has continuum many ergodic measures.
For $A=\{0,1,3\}$ the mask $1+z+z^3$ differs from its
reciprocal by $z(1-z)$ and has no unit-circle zero.  Hence
$X_A=\{x_*\}$ with $x_*\equiv1/3$.  These two examples realize both
possible counts in~\cref{cor:ergodic-measure-count} when $X_A^{(1)}=\varnothing$.

\section*{AI disclosure}

The authors used LLM tools to search the literature,
check grammar, and improve the writing of this manuscript.  Discussions
with LLM tools helped the authors discover the matrix
identity~\eqref{eq:master-identity} and find the example of two
topologically conjugate shifts in Section~\ref{sec:normalization-conjugacy}.
The authors take full
responsibility for the mathematical arguments and the content of this paper.


\begin{thebibliography}{99}

\bibitem{AkhmedovFulghesu}
A. Akhmedov and D. Fulghesu,
\emph{Arithmetic sets in groups},
Math. Z. \textbf{292} (2019), no.~3--4, 1195--1206.
\href{https://doi.org/10.1007/s00209-018-2125-y}{doi:10.1007/s00209-018-2125-y}.

\bibitem{ANNRW}
A. Andoni, A. Naor, A. Nikolov, I. Razenshteyn, and E. Waingarten,
\emph{H\"older homeomorphisms and approximate nearest neighbors},
in Proceedings of the 59th Annual IEEE Symposium on Foundations of
Computer Science (FOCS 2018), IEEE Computer Society, 2018, 159--169.
\href{https://doi.org/10.1109/FOCS.2018.00024}{doi:10.1109/FOCS.2018.00024}.

\bibitem{BhatiaParthasarathy}
R. Bhatia and K. R. Parthasarathy,
\emph{Positive definite functions and operator inequalities},
Bull. London Math. Soc. \textbf{32} (2000), no.~2, 214--228.
\href{https://doi.org/10.1112/S0024609399006797}{doi:10.1112/S0024609399006797}.

\bibitem{Bhattacharya}
S. Bhattacharya,
\emph{Periodicity and decidability of tilings of $\Z^2$},
Amer. J. Math. \textbf{142} (2020), no.~1, 255--266.
\href{https://doi.org/10.1353/ajm.2020.0006}{doi:10.1353/ajm.2020.0006}.

\bibitem{FuSong}
X.-Y. Fu and Z.-J. Song,
\emph{Rational spectra for finite Hadamard pairs and bounded spectral
sets on the real line},
arXiv:2607.19858v1, 2026.
\href{https://arxiv.org/abs/2607.19858v1}{arXiv:2607.19858v1}.

\bibitem{GreenfeldTao}
R. Greenfeld and T. Tao,
\emph{A counterexample to the periodic tiling conjecture},
Ann. of Math. (2) \textbf{200} (2024), no.~1, 301--363.
\href{https://doi.org/10.4007/annals.2024.200.1.5}{doi:10.4007/annals.2024.200.1.5}.

\bibitem{HermleKreidler}
P. Hermle and H. Kreidler,
\emph{A Halmos--von Neumann theorem for actions of general groups},
Appl. Categ. Structures \textbf{31} (2023), Article 38.
\href{https://doi.org/10.1007/s10485-023-09743-y}{doi:10.1007/s10485-023-09743-y}.

\bibitem{KMMS}
G. Kiss, D. Matolcsi, M. Matolcsi, and G. Somlai,
\emph{Tiling and weak tiling in $(\Z_p)^d$},
Sampl. Theory Signal Process. Data Anal. \textbf{22} (2024), Article 1.
\href{https://doi.org/10.1007/s43670-023-00073-7}{doi:10.1007/s43670-023-00073-7}.

\bibitem{KLMS}
G. Kiss, I. Londner, M. Matolcsi, and G. Somlai,
\emph{A lonely weak tile},
Expo. Math. \textbf{44} (2026), no.~1, 125636.
\href{https://doi.org/10.1016/j.exmath.2024.125636}{doi:10.1016/j.exmath.2024.125636}.

\bibitem{KLM}
M. N. Kolountzakis, N. Lev, and M. Matolcsi,
\emph{Spectral sets and weak tiling},
Sampl. Theory Signal Process. Data Anal. \textbf{21} (2023), Article 31.
\href{https://doi.org/10.1007/s43670-023-00070-w}{doi:10.1007/s43670-023-00070-w}.

\bibitem{KolountzakisSurvey}
M. N. Kolountzakis,
\emph{Orthogonal Fourier analysis on domains},
Expo. Math. \textbf{44} (2026), Article 125629.
\href{https://doi.org/10.1016/j.exmath.2024.125629}{doi:10.1016/j.exmath.2024.125629}.

\bibitem{KonyaginLaba}
S. Konyagin and I. \L aba,
\emph{Spectra of certain types of polynomials and tiling of integers
with translates of finite sets},
J. Number Theory \textbf{103} (2003), no.~2, 267--280.
\href{https://doi.org/10.1016/j.jnt.2003.06.006}{doi:10.1016/j.jnt.2003.06.006}.

\bibitem{Lai}
C.-K. Lai, \emph{Private communication}.

\bibitem{Lang}
S. Lang,
\emph{Algebra}, revised 3rd ed.,
Graduate Texts in Mathematics 211,
Springer, New York, 2002.
\href{https://doi.org/10.1007/978-1-4613-0041-0}{doi:10.1007/978-1-4613-0041-0}.

\bibitem{LevMatolcsi}
N. Lev and M. Matolcsi,
\emph{The Fuglede conjecture for convex domains is true in all dimensions},
Acta Math. \textbf{228} (2022), no.~2, 385--420.
\href{https://doi.org/10.4310/ACTA.2022.v228.n2.a3}{doi:10.4310/ACTA.2022.v228.n2.a3}.

\bibitem{Newman}
D. J. Newman,
\emph{Tesselation of integers},
J. Number Theory \textbf{9} (1977), no.~1, 107--111.
\href{https://doi.org/10.1016/0022-314X(77)90054-3}{doi:10.1016/0022-314X(77)90054-3}.

\bibitem{Rudin}
W. Rudin,
\emph{Fourier Analysis on Groups},
Interscience Tracts in Pure and Applied Mathematics 12,
Interscience Publishers, New York--London, 1962.

\bibitem{Schur}
J. Schur,
\emph{Bemerkungen zur Theorie der beschr\"ankten Bilinearformen mit
unendlich vielen Ver\"anderlichen},
J. Reine Angew. Math. \textbf{140} (1911), 1--28.
\href{https://doi.org/10.1515/crll.1911.140.1}{doi:10.1515/crll.1911.140.1}.

\bibitem{Widder}
D. V. Widder,
\emph{Functions harmonic in a strip},
Proc. Amer. Math. Soc. \textbf{12} (1961), no.~1, 67--72.
\href{https://doi.org/10.1090/S0002-9939-1961-0132838-8}{doi:10.1090/S0002-9939-1961-0132838-8}.

\end{thebibliography}
\end{document}